\date{}
\title{On the scaling limit of planar
self-avoiding walk}
\author{Gregory F. Lawler\footnote {Duke University and
Cornell University;  %
partially supported by the  %
National Science Foundation and the Mittag-Leffler  %
Institute.}\and Oded Schramm\footnote {Microsoft Research.}\and  %
Wendelin Werner\footnote{Universit\'e Paris-Sud and IUF.}}  %
\documentclass[12pt]{article}
\usepackage{amsmath}
\usepackage{amssymb}
\usepackage{amsthm}
\usepackage{amsfonts}
\usepackage{graphicx}
\usepackage{amsmath}
\usepackage{amsthm}
\usepackage{amsfonts}
\newif\ifhyper\IfFileExists{hyperref.sty}{\hypertrue}{\hyperfalse} 
\ifhyper\usepackage{hyperref}\fi 
\newtheorem{theorem}{Theorem}
\newtheorem{conclusion}{Prediction}
\newtheorem{corollary}{Corollary}
\newtheorem{lemma}[corollary]{Lemma}
\newtheorem{proposition}[corollary]{Proposition}

\newcommand{\Prob} {{\bf P}}
\newcommand{\Z}{{\mathbb Z}}

\newcommand{\R}{{\mathbb R}}
\newcommand{\C}{{\mathbb C}}

\newcommand{\dist}{{\rm dist}}

\def \eps {\varepsilon}

\def \P {{\bf P}}

\def \H {{\mathbb H}}

\def \Disk {{\mathbb U}}
\def \diam {{\rm diam}}
\def \p  {\partial}

\def \musaw {{\mu_{SAW}}}
\def  \musap {{\mu_{SAP}}}
\def \hmusap {{\hat \mu_{SAP}}}
\def \msaw {{m_{SAW}}}
\def \msap {{m_{SAP}}}
\def \mhsap{{\hat m_{SAP}}}
\def \mbarsap {{\bar m_{SAP}}}
\def \rest{{\mathrm{\SCR/}}}
\def\SCR/{CR}

\def \br {{\rm rloop}}
\def  \ubr {{\rm uloop}}

\def\bloop{\mathrm{bloop}}
\def \U {{\Disk}}
\def \bre {{\medbreak \noindent}}

\begin{document}
\maketitle
 
\begin{abstract}
A planar self-avoiding walk (SAW) is a nearest neighbor random walk
path in the square lattice
with no self-intersection.  A planar self-avoiding
polygon (SAP) is a loop with no self-intersection.
In this paper we present conjectures for the 
scaling limit of the uniform measures on these objects.
The conjectures are based on recent results 
on the stochastic Loewner evolution and non-disconnecting
Brownian motions.  New
heuristic derivations are given for the critical
exponents for SAWs and SAPs.
\end{abstract}

\section{Introduction}

A self-avoiding walk (SAW) in the integer lattice $\Z^2$
starting at the origin of length $n$ is a nearest neighbor
path $\omega = [\omega_0 = 0,\omega_1,\ldots,\omega_n]$ that
visits no point more than once.  Self-avoiding walks were
first introduced as a  
lattice model for polymer chains.  While the model is
very easy to define, it is extremely difficult to 
analyze rigorously. 

Let $C_n$ denote the  number of SAWs of length $n$
starting at the origin. It is easy to see that
$C_{n+m} \leq C_n\: C_m $, and hence subadditivity
arguments imply that there is a $\beta$, called
the {\em connective constant}, such that
\[          C_n \approx \beta^n , \;\;\;\;\;
                     n \rightarrow \infty.   \]
Here, and throughout this paper, we use $\approx$
to mean that the logarithms are asymptotic, so this
can be written 
$ \lim_{n \rightarrow \infty}  n^{-1}  {\log C_n}
              = \beta$. 
The value of $\beta$ is not known exactly and 
perhaps never will.  It  is known rigorously
to be between $2.6$ and $2.7$ and has been estimated
nonrigorously as $\beta =2.638\cdots$ (see
\cite[Section 1.2]{MS}).   The connective constant
is a {\em nonuniversal} quantity in that it varies from
lattice to lattice; it is believed that 
$\beta = \sqrt{2 + \sqrt 2}$ on the honeycomb lattice
\cite{N1}.
It is conjectured that
\[             C_n \sim c\; n^{\gamma -1} \; \beta^n,
      \;\;\;\;\;  n \rightarrow \infty. \]
Here $\sim$ denotes asymptotic, $c$ is a constant, and
$\gamma$ is a nontrivial {\em critical exponent}.
This critical exponent (as well as other critical exponents
for SAWs) is expected to be a {\em universal} quantity;
in particular, unlike the connective constant, the 
exponent will not change if we consider SAWs on
a different (two-dimensional) lattice.  Nienhuis \cite{N1}
gave the exact (but nonrigorous) prediction $\gamma =
43/32$ using a renormalization group technique. This
value was rederived (again nonrigorously) using conformal
field theory ideas  (see \cite{D1,D2}).  Another
critical exponent $\nu$ is defined by the
relation
\begin {equation}
\label {nudef}
           \langle |\omega(n)|^2 \rangle 
               \approx n^{2 \nu}, \;\;\;\;
            n \rightarrow \infty . 
\end {equation}
Here $\langle \cdot \rangle$ denotes the expected value
over the uniform measure on SAWs of length n.  Flory
\cite{Flory} was the first to conjecture $\nu = 3/4$;
the nonrigorous renormalization group and conformal field
theory approaches also give this exact value. None
of these approaches give a prediction for the scaling
limit of SAWs as a measure on continuous paths.

Mandelbrot \cite{mandel}
considered simulations of planar simple random walk
loops (simple random walks conditioned to begin and
end at the same point) and looked at the ``island''
or ``hull'' formed by filling in all of the 
area that is disconnected from infinity.
 He observed that the boundary of these hulls
            were simple curves with fractal dimension
            about $ 4/3$.  Since the conjecture $\nu = 3/4$ suggests
            that paths in a scaling limit of SAWs should have
            dimension $4/3$, he proposed the outer boundary
            of a Brownian loop as a possible candidate for
            what he called``self-avoiding Brownian motion.''
            In particular, he conjectured that the Hausdorff
            dimension of the outer boundary of planar Brownian
            motion is $4/3$

Recently, the authors \cite{LSW1,LSW2,LSWanal} 
proved Mandelbrot's conjecture
for the dimension of the planar Brownian frontier.  
The proof makes use of a new process introduced
by Schramm \cite{S1} called the {\em stochastic
Loewner evolution} ($SLE_\kappa$).  This is a one
parameter family of conformally invariant processes,
indexed by $\kappa$.  It turns out that conformal
invariance of the Brownian hull combined with 
a certain ``locality'' property determine its
distribution.  This distribution is essentially
the same as that of
the hull derived from $SLE_6$.  The hull of
$SLE_6$ can be analyzed by studying a relatively
simple one
dimensional stochastic differential equation and
the dimension of the frontier can be derived
from this.

The main purpose of this paper is to study the 
current status of results relating SAWs and SLEs
(and to planar Brownian frontiers).
In particular:

\begin {itemize}
\item
We give precise formulations  of the existence 
of the scaling limit conjecture, i.e.,
we explain which measures on 
discrete SAWs should converge in the scaling limit.
\item
 We point out that {\em if
the scaling limit of SAWs exists and is conformally
covariant}, then the scaling limit of SAWs is
 $SLE_{8/3}$. 
\item
 Using rigorous exponent calculations
for $SLE_{8/3}$, we give (nonrigorous)
predictions for the exponents $\gamma,\nu$.  In this
way, we rederive the predictions $\nu = 3/4, \gamma = 43/32$.
At the present moment, we do not know how to prove the existence of
the limit (nor, of course, how to prove
that it is conformally covariant).  

\item
We also give precise forms to the conjectures of 
Mandelbrot that the scaling limit of SAPs is the
outer boundary of loops of Brownian motion.  

\end {itemize}

Recently,   Tom Kennedy
\cite{Ken} did Monte Carlo tests of our prediction that
the scaling limit of SAW is $SLE_{8/3}$ and found
excellent agreement between the two distributions.

SAWs have been considered
on the integer lattice $\Z^d$ and the behavior depends
on the dimension.  For $d > 4$, Hara and Slade proved
that the scaling limit of  (the uniform measure on) SAWs
is the same as the scaling limit of simple random walk,
i.e., 
  Brownian motion.
For $d = 4$, it is conjectured that the limit is
also Brownian motion, with a logarithmic scaling
correction, and in $d=3$, the SAW is expected to have
different exponents.  See \cite{MS} for results on
SAWs.  Since we are interested in properties of
SAWs that arise from conformal invariance, we restrict
this paper to $d=2$ (except that the result in the
appendix   is valid for all $d$).

We divide this paper into three parts.  In the next section,
we will discuss rigorous results about $SLE_{8/3}$.  Proofs
of most important facts will appear in \cite{LSWrest}. 
In Section \ref{sawsec}, we define and discuss conjectures
about critical exponents and scaling limits for SAWs and
SAPs.  Here we are not rigorous, but we are
  precise in formulating conjectures, in particular,
on the assumption of conformal invariance or conformal
covariance, the latter term being used for measures
that are invariant up to a scaling factor.  Sections
\ref{slesec2} and \ref{sawsec} are independent of each
other.  In  Section~\ref{slesawsec}, we use the rigorous
properties of $SLE_{8/3}$ and the conjectured properties
of the scaling limit of SAWs to conjecture what the
scaling limit must be.  In particular, we conclude
(nonrigorously) that the continumm limit of SAWs must
be  $SLE_{8/3}$. We then use the rigorous values of
exponents for $SLE_{8/3}$ to
derive  (nonrigorously) the exact values for 
  $\nu,\gamma$, and a third ``half plane''
exponent $\rho$.   In that section we label
a number of statements as ``Predictions''.  These
are not to be interpreted as rigorous theorems but rather
as deductions from nonrigorous arguments.  In other words,
they are conjectures  with significant heuristic
justification.  Finally, the appendix adapts an argument of
Madras and Slade~\cite{MS} to prove that the large $n$ limit
of half-space SAW exists.

\medskip \noindent{\bf Acknowlegment.}
We are grateful to Neal Madras for
helpful advice regarding the existence of the large $n$ limit of the
half plane self avoiding walk.  

\section{SLE and restriction}  \label{slesec2}

\subsection{Motivation}

One of the physicist's key starting points for
studying the behavior of long SAW is the asumption (or conjecture) 
that it has a conformally invariant scaling limit. 
Bearing in mind what conformal invariance means for 
planar Brownian motion (or for 
those models for which it has now 
been proved to hold \cite {Smirnov, LSWlesl}),
 one sees that it is actually more natural 
to forget about the particular
 time parametrization and to study 
paths that are stopped at an exit time (and possibly conditioned 
on the value of their endpoint). 
For instance, the scaling limit of random walk 
 started at a given point and 
stopped at the first exit time of a domain is stopped Brownian motion.
Furthermore, the image of this path  
under a conformal map is another stopped
Brownian moton (more precisely, the image can
be given a parameterization so that it is a
Brownian motion, see, e.g. \cite[Theorem V.1.1]{Bass}).  

It can heuristically be argued (see \S\ref {secheuristic}),
that scaling
limits   of long self-avoiding walks 
remain self-avoiding. 
It is therefore natural  to study measures on continuous 
self-avoiding paths with prescribed endpoints (but non-prescribed
length).

Suppose $D$ is a bounded simply connected domain
in $\C$ and $z,w \in \partial D$. 
 Assume that
$\partial D$ is smooth near $z$ and $w$.  Our goal
is to find a measure that will be a candidate for
the scaling limit of a measure on SAWs from $z$ to $w$
in a grid approximation of $D$.
This should be a measure on simple curves from $z$ to
$w$ staying in $D$ that is conformally invariant
and also (in some sense) a ``uniform measure'' on
such curves.

A notion of ``uniform measure''  requires
some reference measure 
and at this point we do not
know what this reference measure would be, so this 
is rather vague.  However, a uniform measure 
has the property that if it is restricted to a smaller set, it remains
uniform. 
For each $D, z, w$ as above, let $m^\# (z,w;D)$ denote such a 
``uniform'' measure (we will use the superscript $\#$ to indicate
that these are 
probability measures).   
Then if $D' \subset
D$ and is smooth near $z$ and $w$,
the restriction of $m^\#(z,w;D)$ to
those curves that lie in $D'$ is also ``uniform''.
Hence the conditional law of $\gamma$ restricted to $\{ \gamma \subset D'\}$
should be  $m^\# (z,w; D')$.
We say that a family of probability measures that satisfy this
property is {\em restriction covariant}.

The conformal invariance asumption is naturally 
phrased as follows: 
For all $(z,w;D), (z', w'; D')$ and all conformal map 
$f$ from $D$ onto $D'$ with $f (w) = w' $ 
and $f (z) = z'$, the image of the measure 
$m^\# (z,w;D)$ under $f$ (that we will
denote by  $f \circ m^\#(z,w;D)  
$) should be equal to 
$m^\#(z',w'; D')$. 
We then say that such a  family of measures is 
{\em conformally invariant}
(we use the term covariant when there is a 
multiplicative term, and invariance when there
is identity between the measures).

Note that such a conformally invariant family of 
probability measures is completely determined by  
 $m^\# (0,\infty; \H)$ (for instance) where
$\H$ denotes the upper half-plane
since all other $m^\# (z,w;D)$ are just the
conformal images of this measure under a mapping $\Phi$ from
$\H$ onto $D$ with $\Phi(0)=z$ and $\Phi(\infty)=w$.
(Note that a hull connecting $0$ and $\infty$ in $\H$
is unbounded, but we think of it as a compact subset of
$\C\cup\{\infty\}$.)

As we shall see in the next subsection, there exists 
in fact only one family of probability measures $m^\# (z,w;D)$ 
on simple curves that is conformally invariant and 
restriction covariant.

\subsection{Restriction measures (chordal)}

In order to find such measures on self-avoiding curves,
it is convenient to work in a slightly 
more general setup, and to consider measures on 
{\em hulls.}
If $D$ is a simply connected domain and $z,w \in
\p D$ (such that $\p D$ is smooth near $z$ and $w$),
 then a {\em hull} connecting $z$ and $w$
in $D$ is a compact set $K$ with $K \subset
D \cup \{z,w\}$  such that $D\setminus K$ has
exactly two connected components.  Note that the
image of a simple curve in $D$ connecting $z$ and $w$ is
a hull.  Moreover, if a finite collection of
(not necessarily simple) curves connecting $z$ and
$w$ in $D$ is given, then the curves generate
 a hull  $K$ by ``filling
in the holes'' of the union of
the paths (more precisely, $D \setminus K$ consists
of the two connected components of $D$ minus the 
curves whose boundaries include nontrivial arcs in
$\p D$).

Let ${\cal D}$ denote the class of all triplets $(D,z,w)$ such that 
$D$ is bounded simply connected, $z \not=w \in \partial D$, and $\partial D$
is smooth in the neighborhood of $z,w$.

\medskip\medskip

\noindent {\bf Definitions}.  Suppose that
for each $(D,z,w) \in {\cal D}$, we have a probability measure $m^\#(z,w;D)$ on
 hulls from $z$ to $w$ staying in $D$.  We call this
family of measures a {\em conformal  restriction (\SCR/) family}
if it satisfies:
\begin{itemize}
\item Restriction covariance:
if $D' \subset D$,
$z,w \in \p D \cap \p D'$ and $\p D,\p D'$ agree
near $z,w$, then $m^\# (z,w;D)$ restricted to 
the set of hulls
that lie in $D'$ is a multiple of the measure $m^\# (z, w;D')$.

\item Conformal invariance: if   
$f: D \rightarrow D'$ is a conformal
transformation,
\[         f \circ m^\#(z,w;D) = m^\#(f(z),f(w);D') . \]
\end{itemize}

\medskip\medskip

Suppose $m^\#(z,w;D)$ is a \SCR/
family.  By conformal invariance, the family is
determined by $m^\#(0,\infty;\H)$, where $\H$ denotes
the upper half plane. Conformal invariance 
implies that $m^\#(0,\infty;\H)$ is
invariant under the maps $z \mapsto rz (r > 0)$.
Suppose $A$ is a compact set not containing
the  origin
such that $\H \setminus A$ is simply connected.  
 By considering infinitesimal $A$,
and using conformal invariance and the restriction
property, the following can be proved.

\begin{lemma}[\cite{LSWrest}]  \label{feb7prop}
If $m^\#(z,w;D)$ is a \SCR/ family, then there
is an $a > 0$ such that the $m^\#(0,\infty;\H)$-probability
 that a hull lies in $\H \setminus
A$ is $\Phi'_A(0)^a$, where 
$\Phi_A$ is the unique conformal transformation
of $\H \setminus A$ onto $\H$ with $\Phi_A(0) = 0,
\Phi_A(\infty) = \infty, \Phi_A(z')  \sim z' $ when $z' \to \infty$.
\end{lemma}

Conversely, suppose that
$m^\#(0,\infty;\H)$ is a probability measure on
hulls such that
for all $A$ as above, the probability that a hull lies
in $\H \setminus A$ is $\Phi_A'(0)^a$.
This uniquely determines $m^\#(0,\infty;\H)$ \cite{LSWrest}.
Define probability measures $m^\#(z,w;D)$ by
conformal transformation.  Then $m^\#(z,w;D)$
gives a \SCR/ family that denote by
$\rest(a)$.  The measure $m^\#(z,w;D)$
will then be denoted by $\rest^a(z,w;D)$.

It follows from the lemma that every \SCR/ family satisfies
the symmetry relation $\rest^a(z,w;D)=\rest^a(w,z;D)$ and
that $\rest^a(z,w;D)$ is invariant under the anti-conformal
automorphisms of $D$ which leave $z$ and $w$ invariant.

One can also define the corresponding
measures 
$$m(z,w;D) := |f'(z)|^a \: |f'(w)|^a
\; m^\#(z,w;D),$$
 where $f$ is a conformal
transformation of $D$ onto the unit
disk $\Disk$ with $f(z) =1,
f(w) = -1$. 
Although the map $f$ is not
uniquely defined, the quantity $|f'(z)| \: |f'(w)|$
does not depend on the choice of $f$. The
family $m(z,w;D)$ then satisfies the following properties:
\begin{itemize}
\item 
Restriction invariance:
 if $D' \subset D$,
$z,w \in \p D \cap \p D'$, $\p D,\p D'$ agree
near $z,w$,  then $m(z,w;D')$ is equal to $m(z,w;D)$
restricted to hulls that stay in $D'$.
\item Conformal covariance:
if   
$f: D \rightarrow D'$ is a conformal
transformation,
then $f \circ m(z,w;D) $ is a multiple of $m (f(z), f(w); D')$.
In fact,
\[ f \circ m (z,w;D) 
= |f'(z)|^a \: |f'(w)|^a\;
m(f(z),f(w);D') . \]
\item Symmetry: $m(z,w;D) = m(w,z;D)$.
\end{itemize}

Using Brownian motion, 
it is not difficult \cite {LSWrest} (see also~\cite{LW1})
to construct the $\rest (k)$ 
families for integer $k$: 
\begin{itemize}
\item  The $\rest(1)$ family is obtained from a 
 half plane Brownian
   excursion, i.e., a Brownian motion started at the origin
   conditioned never to leave the upper half plane (this can
    be obtained by letting the $x$-component move according
    to a Brownian motion and the $y$-component according
     to an independent
    Bessel-3 process, which is the same as the $l^2$ norm
    of a three dimensional Brownian motion).
\item  If $m^\#(z,w;D)$ is an $\rest(a)$ family and
     $m_1^\#(z,w,D)$ is an $\rest (a_1)$ family then
       $m^\#(z,w;D) \times m_1^\#(z,w,D)$ gives an
      $\rest(a + a_1)$ family.  Here, the hull is 
        obtained from the union of the two independent
        hulls by ``filling in the holes''.  In particular, for
       positive integers $k$, the $\rest(k)$ family is obtained
       from the hull of  the union of $k$ independent Brownian excursions.
\end{itemize}
Note that the $\rest(k)$ families for integer $k$ are not supported on
simple paths.
In \cite {LSWrest},
all $\rest(a)$ families are determined
and studied. It turns out that 
  there is only one value of
$a$ such that $\rest(a)$ is supported on simple paths.

\begin{theorem}[\cite{LSWrest}] 
\label {th1}
(i)
The $\rest(a)$
family  exists for
all $a \geq 5/8$ and for no $a < 5/8$.  

(ii)
$\rest(5/8)$ is the only \SCR/ family supported on simple paths.
The measure $\rest^{5/8}(0,\infty;\H)$ is 
chordal $SLE_{8/3}$.
\end{theorem}

We describe chordal $SLE_{8/3}$ in the next subsection.
 
\subsection{Stochastic Loewner evolution }\label{chordslesec}

In this section, we briefly review the definition
of the stochastic Loewner evolution;
 see \cite{LSW1,LSW2,S1} for more details.
First, let us loosely describe the motivation behind
the definition.  Suppose that we wish to give a 
probability measure on curves $\gamma:[0,\infty)
\rightarrow \overline \H$.  We are interested only
in the curve and not in its parametrization; we
consider two curves to be the same if one is just a
reparametrization of the other. 
Suppose
that an initial piece of the curve $\gamma[0,t]$
is given.  If $\gamma$ is SAW, for example, then the rest of the
path is still SAW, in the domain with $\gamma[0,t]$ removed.
When we combine this Markovian type property with a conformal invariance
assumption, we get the following.
  Suppose $\phi$ is a conformal
transformation of    the
unbounded component of $\H \setminus \gamma[0,t]$ onto
$\H$ with $\phi(\infty) = \infty$ and $\phi(\gamma(t))
= 0$. Then the conditional distribution of
$\phi \circ \gamma(s), t \leq s < \infty$ given
$\gamma[0,t]$ is the same as (a time change of)
the distribution of $\gamma[0,\infty)$.  If we assume
this property and also a symmetry property (the 
distribution is invariant under reflections about
the imaginary axis), then the only possible distribution
for $\gamma$ is chordal
$SLE_\kappa$ as defined below.
 
   Suppose $\gamma:[0,\infty)
\rightarrow \overline \H$ is a curve, and let $\H_t$
be the unbounded component of
$\H \setminus \gamma[0,t]$.  Let $g_t:\H_t
\rightarrow \H$ be the unique conformal transformation
such that $g_t(z) - z = o(1)$ as $z \rightarrow \infty$;
$g_t$ has the expansion
\[   g_t(z) = z + \frac{a(t)}{z} + \cdots , \;\;\;
   z \rightarrow \infty.\]
 We suppose that
$\gamma(z)$ either is simple or ``bounces off'' itself;
more precisely, we assume that $W_t = 
g_t(\gamma(t))$ is a continuous function from
$[0,\infty)$ to $\R$ and the domains $\H_t$ are
strictly decreasing in $t$.
In this case, $a(t)$ is continuous, strictly increasing.
   We can  assume that
$a(t) = 2t$; otherwise, we can reparametrize the
curve so this is true.  Then $g_t$ satisfies
the Loewner differential equation
\[   \p_t g_t(z) = \frac{2}{g_t(z) - W_t} ,
\;\;\;\; g_0(z) = z . \]
Conversely, we can start with a function $W_t$
and obtain the maps $g_t$ by solving the equation
above.  If things are sufficiently nice, then we
can also define $\gamma(t) = g_t^{-1}(W_t)$.
The {\em (chordal) stochastic Loewner evolution with parameter
$\kappa \geq 0$} or {\em  chordal 
$SLE_\kappa$} is the process obtained by choosing
$W_t = B_{\kappa t}$ where $B$ is a standard one-dimensional
 Brownian motion.  It has been shown
\cite{RS,LSWlesl} that the curve $\gamma(t) = g_t^{-1}(W_t)$
is almost surely  well defined  for any   $\kappa$.
  This defines $SLE_\kappa$ on the
upper half plane. It can be defined for other domains
by conformal transformation (this uses the fact
that $SLE_\kappa$ in the upper half plane
is invariant under scaling
in an appropriate sense).

The qualitative behavior of the paths $\gamma$ vary
significantly as $\kappa$ varies (see \cite{RS}).  In particular, for
$\kappa \leq 4$, the $SLE_\kappa$ paths are simple
paths while for $\kappa > 4$ they have self-intersections.
For $\kappa \geq 8$ the curves are space filling.

Let us now list some properties that are specific
to the very special value $\kappa= 8/3$.

\begin{proposition}[\cite{LSWrest}]  \label{restprop}
Suppose $\gamma = \gamma[0,\infty)$ is the path of chordal
$SLE_{8/3}$.

(i)
If $A\subset\C\setminus\{0\}$ is a compact set such
that $\H \setminus A$ is simply connected and
$\Phi_A$ is as in Proposition \ref{feb7prop}, then
the probability that $\gamma(0,\infty) $ lies
in $\H \setminus A$ is $\Phi_A'(0)^{5/8}$.

(ii)
The Hausdorff dimension of $\gamma$ is $4/3$.

(iii)
The hull of the union of eight independent copies
of $\gamma$ has the same law as the law 
of the hull of the union of five Brownian excursions 
in the upper half plane (as defined in the 
previous subsection).

\end {proposition}
  
(i) shows that $SLE_{8/3}$ indeed  gives the  restriction
family with $a=5/8$.  It is also proved in \cite{LSWrest}
that for any other $\kappa \leq 4$, the probability
of staying in $\H \setminus A$ is not given by
$\Phi_A'(0)^\alpha$ for some $\alpha$.  However, these
other values of $\kappa$ are used to construct
the $\rest(a)$ families for $a > 5/8$ and this
construction shows that these families are not
supported on simple paths.    This characterization of
$SLE_{8/3}$ as the only $\rest(a)$
family supported on simple paths will tell us that it is the only possible
scaling limit of SAW {\em assuming the scaling
limit exists and is conformally invariant}

(iii) follows simply from the fact that both hulls give
a realization of $\rest (5)$, which is unique. (ii) then 
follows from the fact  
\cite{LSW1,LSW2,LSWanal}  that the boundary
of the frontier,
or outer boundary,  of planar Brownian motion has
dimension $4/3$. 

(iii) shows that the frontier of a Brownian motion look
locally the same as an $SLE_{8/3}$ and hence --- provided 
that $SLE_{8/3}$ is in some sense the scaling limit of SAW --- 
this agrees with Mandelbrot's visual observation.

\subsection {Some consequences} \label{conseqsec}

The relation given by
Proposition \ref {restprop}(iii)  
between Brownian excursions and $SLE_{8/3}$ makes it possible 
to show that some other measures defined in terms of $SLE_{8/3}$ or Brownian
excursions are equal.

\bre
{\bf Boundary loops.}
Motivated by possible scaling limits of self-avoiding walks
in a domain that start and end at the same boundary point,
one can try to define an $SLE_{8/3}$ in $D$ that starts
 and ends at $z$. It turns 
out that the natural measure is in fact an infinite measure that one can view
as the appopriately rescaled limit of $\rest^{5/8}(z,w;D)$ 
as $w \to z$.
In the upper half-plane, one defines
$$
\mu_{\bloop} (0; \H) = \lim_{\eps \to 0} \eps^{-2}\rest^{5/8}(0,\eps;\H).
$$
The relation between $SLE_{8/3}$ and the Brownian excursion shows easily that if 
one considers the limit when $\eps$ goes to $0$ of $\eps^{-2}$ times
the law of the ``outer frontier'' of a Brownian excursion in $\H$
from $0$ to $\eps$, then one gets exactly $\frac 85\, \mu_{\bloop}$.
We call $\mu_{\bloop}$ the {\em boundary loop measure}.

Another way to describe the measure on boundary
loops at zero in the upper half plane
is to give the measure of the set of hulls
that hit a set $A$.  Suppose $A\subset\C\setminus\{0\}$ is compact and
$\H \setminus A$ is simply connected.  Let
$\Phi_A$ be the map as in Lemma~\ref{feb7prop}.
Then the measure of the set of bubbles that hit $A$
is 
$$
-\frac 5{48} \,\mathcal S_{\Phi_A}(0)\,,
$$
where
$$
\mathcal S_f:=\frac {f'''}{f'}-\frac {3{f''}^2}{2{f'}^2}
\,,
$$ is known as the
{\em Schwarzian derivative} of $f$.

\bre
{\bf Non-intersecting $SLE_{8/3}$'s.}
Similarly, one can consider the law of the hull of 
two $SLE_{8/3}$ that are ``conditioned not to intersect''.
The procedure here is to take the law of two 
(independent) $SLE_{8/3}$'s $\gamma$ and 
$\gamma'$, one from $z$ to $w$ in $D$ and
the other one from $z'$ to $w'$ in the same domain $D$ (where
$z, z', w', w$ are ordered clockwise on $\partial D$ say)
and to condition these paths by $\gamma \cap \gamma' = \emptyset$.
Then, one looks at the limit of this measure when $z' \to z$ and $w' \to w$.
It turns out that the obtained law is exactly the $\rest (2)$ 
family, i.e., it is the same as considering the union of two
independent Brownian excursions from $z$ to $w$.

\subsection{Radial case}  \label{radslesec}
\label {S.radial}
One can raise analogous questions considering random paths
that join an interior point of a domain to a 
boundary point that should appear as scaling limits of 
some SAW.
We will not do the whole analysis of restriction measures
in this setup here, but only briefly state some results
concerning $SLE_{8/3}$ and its 
relations to Brownian frontiers. For symmetry reasons,
the most natural domain to work in 
is the unit disk. One wishes to 
study  random paths from $0$ to $1$ in the
unit disk.

\bre
{\bf Radial $SLE_{8/3}$.}
A version of the stochastic
Loewner evolution called {\em radial SLE}  defines
 such a path. 
 Consider the initial
value problem
\[    \p_t \hat g_t(z) = \hat
g_t(z) \: \frac{e^{i W_t} + \hat g_t(z)}
                 {e^{i W_t} - \hat g_t(z)} , 
\;\;\;\; \hat g_0(z) = z\]
where $W_t = B_{\kappa t}$ as above.  If we write
$\hat g_t(z) = \exp\{i \hat h_t(z)\}$ (in a neighborhood of
$\hat g_t(z)$) we can write this as
\[   \p_t \hat h_t(z) = \cot[\frac{\hat h_t(z) - W_t}{2}].\]
 If $\kappa \leq 4$,
then
$\hat \gamma(t) := \hat g_t^{-1}(e^{i W_t})$ is a simple curve
going from a
point on the unit circle to the origin staying
in the unit disk $\Disk$, and  $\hat g_t$ maps
 $\Disk \setminus \hat \gamma[0,t]$
  conformally onto $\Disk$ with
$\hat g_t(0) = 0, \tilde g_t'(0) = e^t$.  
There is a close relationship between chordal
$SLE_\kappa$ and radial $SLE_\kappa$ for the same
value of $\kappa$ \cite[Theorem 4.1]{LSW2};
for example, the Hausdorff dimensions
of the paths are the same.  There is a special
relationship that holds for $\kappa = 8/3$: radial
$SLE_{8/3}$ evolves like  the image of chordal $SLE_{8/3}$ under the
map $z \mapsto e^{iz}$ weighted so that the paths never
form a loop about the origin \cite{LSWrest,LSWrest2}.

\begin {theorem}
The radial $SLE_{8/3}$ is the unique stochastic Loewner 
evolution in $\U$ (from $1$ to $0$) such that 
for all
simply
connected subdomain $U$ of $\Disk$ containing
the origin whose boundary contains
an open subarc of $\p \Disk$ containing $1$ the
event $\gamma\subset\U$ has positive probability
(where $\gamma$ denotes the $SLE$ path) and if
$\Psi:U \rightarrow 
\Disk$ is the conformal transformation with
$\Psi(U)=\Disk, \Psi(0) = 0, \Psi(1) = 1$, then the law of $\Psi (\gamma)$
conditionally on $\{ \gamma \subset U\}$
 is equal to 
the law of $\gamma$ itself.  In fact, 
\begin {equation}
\label {548}
\P [ \gamma \subset U ] 
=|\Psi'(1)|^{5/8} \: |\Psi'(0)|^{5/48}.
\end {equation}
\end {theorem}

This is, of course, closely related to restriction covariance and
conformal invariance in the chordal case.
Note the different exponents associated to the interior point 
and to the boundary point.

   Let $m^\#(0,1;
\Disk)$ denote the probability
measure on curves in $\Disk$
  given
by radial $SLE_{8/3}$ started at $1$.  
  By conformal
transformation, we then get a family of
probability measures
$m^\#(w,z;D)$ where $D$ ranges over simply connected
domains, $z$ over points in $\p D$ and $w$ over
points in $D$. As before, this family
is restriction covariant and conformally invariant, with the
appropriate modifications of the definitions to the radial setting.

If $f: \Disk \rightarrow D$
is a conformal transformation with $f(\Disk)=D, f(1) = z,
f(0) = w$ and $f$ is smooth at $1$,
we let
$$m(w,z;D) = |f'(1)|^{5/8} \; |f'(0)|^{5/48}\;
m^\#(w,z;D) = |f'(1)|^{5/8} \; |f'(0)|^{5/48}
f \circ m^\#(0,1;\Disk).$$
Just as in the chordal case, the measures
$m(z,w;D)$  are then conformally covariant and restriction invariant.

Note that we have used the same notation here as for (chordal)
$\rest$
families; however, in this case $w$ is an interior
point so there is no inconsistency in notation.

\bre
{\bf The full-plane $SLE$.}
Using the inversion $z \mapsto 1/z$, we can also define
radial $SLE_{\kappa}$ going from the unit circle to infinity
staying in $\C \setminus \Disk$.  By letting the unit
circle shrink to a point, we get a limiting process which
is called {\em full plane} $SLE_{\kappa}$.  This 
can also be defined directly.  Let $B_t, -\infty < t < \infty$,
be a standard Brownian motion such that $B_0$ has a
uniform distribution on $[0,2 \pi]$, and let
$W_t = B_{\kappa t}$.  Full plane 
$SLE_\kappa$ driven by $e^{i\,W_t}$ is
the unique curve 
$\hat\gamma(t)$, $- \infty < t < \infty$,
with  $\hat \gamma(-\infty) =0,
\hat \gamma(\infty) = \infty$ satisfying the following.
Let $g_t$ be the conformal transformation of
the unbounded component of $\C \setminus 
\hat \gamma[-\infty,t]$
onto $\C \setminus \overline \Disk$ with $ g_t(\infty)=\infty$
and $g_t'(\infty)>0$. Then $g_t'(\infty) = e^{-t}$
and $g_t$ satisfies
\[    \p_t g_t(z) = g_t(z) \; \frac{e^{i\,W_t} + g_t(z)}{ 
   e^{i\,W_t}- g_t(z)} . \]
The transformations $g_t$ can be obtained as
the limit as $K \rightarrow -\infty$
of the solution $g_t^K$, $K\leq t < \infty$, of this
equation with initial condition $g_t^K(z) = e^K \: z.$
The path $\hat\gamma(t)$ can be determined from $W_t$
by $\hat\gamma(t)=g_t^{-1}(e^{iW_t})$.


\bre
{\bf Non-disconnecting Brownian motions.}
A ``typical'' point on the frontier of a planar
Brownian motion looks like a pair of Brownian motions
(the ``past'' and the ``future'')
starting at the point
conditioned so that the union of the paths does not
disconnect the point from infinity.  This is a conditioning
on an event of measure zero, but it can be made precise
by taking an appropriate limit (see, e.g., \cite{LSWanal}).
For example, one can consider pairs of Brownian motions
starting at the origin; let $B_\epsilon$ be the union
of the two paths from the first time they reach the
circle of radius $\epsilon$ to the first time they
reach the circle of radius $1/\epsilon$.  Let
$V_\epsilon$ be the event that the origin is in
the unbounded component of the complement of $B_\epsilon$.
Then \cite{LSW1,LSW2,LSWanal}
as $\epsilon \rightarrow 0+$, $\Prob(V_\epsilon)
\asymp \epsilon^{4/3}$
(where $\asymp$ means that the ratio of both sides is bounded
away from zero and infinity), 
 and there is a limiting
measure on paths that we call a {\em non-disconnecting
pair of Brownian motions} starting at the origin.

The ``hull'' generated by a pair of non-disconnecting
Brownian motions is bounded by a pair of infinite
curves from the origin to $\infty$.  It is known
\cite{BL} that these curves do not intersect.  (If
we took only one Brownian path, conditioned it not
to disconnect the origin, and then considered the
``left'' and ``right'' paths of the frontier, these
paths {\em would} have self-intersections; this is
related to the fact that planar Brownian motions
have cut points.)

If $D$ is a domain, $z \in \partial D$,
$w \in D$, we can similarly define pairs of Brownian
motions starting at $w$, conditioned to leave $D$ at
$z$ and conditioned so that $w$ is not
disconnected from $\partial D$.  Let $m^\#(w,z;D)$
be the probability measure on curves given by this
construction.  Let
\[  m(w,z;D) = |f'(1)|^2 \; |f'(0)|^{2/3} \;
            m^\#(w,z;D) , \]
where $f: \Disk \rightarrow D$ is a conformal
transformation with $f(1) = z, f(0) = w$.  
The scaling exponent $2/3$ is the
$2$-disconnection exponent $\eta$ for Brownian motion
defined by $\Prob(V_\epsilon) \asymp \epsilon
^{2 \eta}$.  The
measures $m(w,z;D)$ can be seen to satisfy
restriction invariance 
and  conformal covariance,
because
if   
$f: D \rightarrow D'$ is a conformal
transformation,
\[         f \circ m(w,z;D) = |f'(z)|^{2} \: |f'(w)|^{2/3}
\:m(f(w),f(z);D') . \]
It is likely that one can again identify this measure on
hulls of non-disconnecting 
Brownian motions to pairs of radial $SLE_{8/3}$'s  
``conditioned not to intersect''.
By zooming in towards the origin, this may lead to a proof
that the inside and the outside of the
corresponding hull have the same law.

\subsection{More loops}

A {\em (rooted)  Brownian loop} (sometimes called
a Brownian bridge) of time duration
$1$ rooted at $0$
is a Brownian motion with $B_0 = 0$ and $B_1 = 0$.  This
involves conditioning on an event of measure zero, but can easily
be made precise in a number of ways.
For example, $B$ can be defined 
as the weak limit as $\epsilon$ tends to zero of
Brownian paths starting at the origin
conditioned on $|B_1|\leq \epsilon$, or, equivalently,
$B_t:= \tilde B_t-t\tilde B_1$, where
$\tilde B_t$ is ordinary Brownian motion started from $0$.
 Let $\mu_{\br,1}$
denote this probability
measure on loops.  Using Brownian scaling
we can also get the probability measure $\mu_{\br,t}$ for
Brownian loops of time duration $t$ ($B_0 = B_t$).  (Here
$\br$ stands for rooted loop.) We define the
rooted
Brownian loop measure rooted at $0$
to be the infinite measure
\[            \mu_{\br}^0 = \int_0^\infty
    t^{-1} \; \mu_{\br,t} \; dt . \]
If $\mu$ and $\mu'$ are measures on collections of paths,
we write
$$
\mu \doteq \mu'
$$
to mean that $\mu$ is the same as $\mu'$ when we forget the parameterization
of the curve;
that is, $\mu$ and $\mu'$ are identical as measures
on equivalence classes paths up to monotone reparameterization.
If $f(z) = rz$ is a dilation, then
\[
  f \circ \mu_{\br}^0  =   \int_0^\infty
           t^{-1} \;  f \circ \mu_{\br,t} \; dt
     \doteq   \int_0^\infty t^{-1} \; \mu_{\br, r^2 t} \;
     dt = \mu_{\br}^0 .
\]
The second equality uses the standard scaling relation
for Brownian motion.  Hence, $\mu_{\br}$ is invariant under dilations
up to $\doteq$ equivalence
(which is the reason for introducing the weighting by $t^{-1}$).

We can define the Brownian loop measure rooted at $z$,
$\mu_{\br}^z$, by translation.  If $D \subset \C$ is a domain,
we let $\mu_{\br,D}^z$ be $\mu_{\br}^z$ restricted to
loops that stay in $D$.  The family of (infinite) measures
 $\mu_{\br,D}^z$  as $D$ ranges over simply connected
domains and $z \in D$ is an example of a family that
satisfies restriction and is conformally
{\em invariant}, i.e., $f \circ \mu_{\br,D}^z
\doteq \mu_{\br,f(D)}^{f(z)}$, if $f$ is conformal in $D$.

Let
\[   \mu_{\br,D} = \int_D \mu_{\br,D}^z \; dA(z) , \]
where $A$ denotes area.  This is a measure on
continuous curves (rooted loops)
$\gamma:[0,T] \rightarrow
D$ with $\gamma(0) = \gamma(T)$ and
$T = T(\gamma)$ is the ``time duration'' of the
curve.   We can also consider this as
a measure on {\em unrooted loops} where we
consider two curves $\gamma,\gamma^1$ the same
if $T(\gamma) = T(\gamma^1)$ and there exists a
$t_0$ such that $\gamma^1(t) = \gamma(t + t_0)$
for all $t$, with addition modulo $T(\gamma)$.
Equivalent curves have the same time duration;
hence we can talk about the time duration
$T(\gamma)$ of an equivalence class of curves.
Let $\mu_{\ubr,D}^z = T^{-1} \; \mu_{\br,D}^z$
(i.e., the measure whose Radon-Nikodym derivative
with respect to $ \mu_{\br,D}^z$ is $T^{-1}$)
considered as a measure on unrooted loops. Let
\[  \mu_{\ubr,D} = \int_D   \mu_{\ubr,D}^z \; dA(z) .\]
   In \cite{LSWrest} it is proved that
the measure $\mu_{\ubr,D}$ is conformally invariant, i.e.,
if $f: D \rightarrow D'$ is a conformal transformation,
then
\[   f \circ \mu_{\ubr,D} \doteq \mu_{\ubr,D'} . \]
For example, if $f(z) = rz$ is a dilation, then
under Brownian scaling times scales by a factor of $r^2$ under
$f$, and
\begin{eqnarray*}
 f \circ \mu_{\ubr,D} & = &  \int_D T(\gamma)^{-1}  \;
   f \circ \mu_{\br,D}^z \; dA(z) \\
  & \doteq & \int_D [T(f \circ \gamma)]^{-1} \; \mu_{\br,f(D)}^{f(z)}  \;
                      r^2 \; dA(z)  \\
  & = & \int_{f(D)} \mu_{\ubr,f(D)}^{w} \; dA(w) =
                 \mu_{\ubr,f(D)}.
\end{eqnarray*}
 It is also immediate from the definition
that $\mu_{\ubr,D}$ satisfies the following restriction
property: if $D' \subset D$, then $\mu_{\ubr,D'}$
is $\mu_{\ubr,D}$ restricted to loops that stay in
$D'$.

An unrooted Brownian loop gives rise to
a Brownian ``hull'' or ``island'' (the hull generated
by  any compact set $K$ is the complement of the
unbounded component of $\C \setminus K$).
The boundaries of these hulls are Jordan curves,
i.e., homeomorphic images of the unit circle.
(Technically we should say that this
holds {\em with probability one}.  We will
omit this phrase when   describing properties
of Brownian paths.)
It is possible to relate these Brownian hulls to the
boundary loop measure, and therefore also to $SLE_{8/3}$.

\section{Self-avoiding walks} \label{sawsec}

\subsection{Definitions}\label{sawdefsec}

A self-avoiding walk (SAW) of length $n$ in $\Z^2$
is a nearest neighbor walk of length $n$ that visits
no point more than once.
More precisely, it is a
finite sequence  $\omega =[\omega_0,
\omega_1,\ldots,\omega_n]$ with
$\omega_j \in \Z^2$;   $|\omega_{j+1}
- \omega_j| =1, j=0\ldots,n-1$; and $\omega_j
\neq \omega_k, 0 \leq j < k \leq n$.  Let $\Lambda_n$
denote the set of all SAWs of length $n$; this is an
infinite set since there are an infinite number of choices
for the initial point $\omega_0$. Let  $\Lambda_n^*$ be
the set of SAWs with $\omega_0 = 0$; $\Lambda
= \cup_{n=0}^\infty \: \Lambda_n; \Lambda^* =
\cup_{n=0}^\infty\; \Lambda^*_n$.  It is easy to see that
$\#(\Lambda_{n+m}^*) \leq \#(\Lambda_n^*)
\: \#(\Lambda_m^*)$, and hence by subadditivity
there exists
a number $\beta$ called the {\em connective
constant} such that 
$  \#(\Lambda_n^*) \approx \beta^n . $
The exact value
of $\beta$ is unknown, and perhaps will never be
know; rigorously it is known to be between
2.6 and 2.7.

 Let $\musaw$ be the measure
on $\Lambda$ that assigns measure $\beta^{-n}$ to
each $\omega \in \Lambda$; $\musaw$ can also be 
considered as a measure on $\Lambda^*, \Lambda_n,
\Lambda_n^*$.  On $\Lambda_n^*$, it is a uniform measure
with total mass $\beta^{-n}\: \#(\Lambda_n^*)$.
(More generally, we might consider  
the measure  which assigns measure 
$e^{-\alpha n}$ to each $\omega \in \Lambda_n$.
When we   vary $\alpha$, these measures on 
$\Lambda^*$ undergo
a  ``phase transition'' 
at $e^{-\alpha} = 1/\beta$  from a
finite to an infinite measure.  The study of systems
near a phase transition goes under the name of
critical phenomena, and many of the critical exponents
discussed below for SAWs have analogues for other
lattice models in statistical physics.)

A self-avoiding polygon (SAP) of length $2n$ in
$\Z^2$ is a finite sequence   $\omega =[\omega_0,
\omega_1,$ $\ldots,\omega_{2n}]$ with $\omega_j \in \Z^2$;
$|\omega_{j+1}
- \omega_j| =1, j=0\ldots,2n-1$; $\omega_0 =
\omega_{2n}$ ; and $\omega_j
\neq \omega_k, 0 \leq j < k \leq 2n-1$.  There is a
one-to-one correspondence between SAPs of length
$2n$ and SAWs of length $2n-1$ whose final point
is distance one from the initial point.
Let $\Gamma_{2n}$ denote the set of SAPs of length
$2n$ and $\Gamma_{2n}^*$ the set of $\omega
\in \Gamma_{2n}$ with $\omega_0 = 0$. Let
$\Gamma = \cup_{n=0}^\infty \Gamma_{2n},
\Gamma^* = \cup_{n=0}^\infty \Gamma_{2n}^*.$  It
is known that $\#(\Gamma_{2n}^*) \approx
\beta^{2n}$, i.e., the connective constant for
SAPs  is the same as for SAWs \cite[Corollary 3.2.5]
{MS}.  Let
$\musap$ be the measure on $\Gamma$ that assigns
measure $\beta^{-2n}$ to each $\omega \in
\Gamma_{2n}$.  

The above definition of a SAP differs from that
in \cite{MS} in that we have specified which point
in a SAP is the initial (and hence terminal) point and we
have specified an orientation.
If $\omega = [\omega_0,\omega_1,\ldots,$ $\omega_{2n-1},
\omega_{2n} = \omega_0] \in \Gamma_{2n}$ and $j \in
{0,\ldots,2n-1}$, let $T^j\omega=[\omega_j,\omega_{j+1},
\ldots,\omega_{2j},\omega_1,\omega_2,$ $\ldots,\omega_{j-1},
\omega_j]$, i.e.,  $T^j\omega$ traverses the same points
as $\omega$ in the same order; the only difference is
the starting point.  Let $\omega^R = [ \omega_{2n},\omega_{2n-1},
\omega_{2n-1},\ldots,\omega_0]$ be the reversal
of $\omega$.  Let $\hat \Gamma_{2n}$ denote the equivalence
classes $[\omega]$ of $\omega \in \Gamma_{2n}$ generated
by the equivalences
$T^j \omega \sim \omega$ and $\omega^R \sim 
\omega$;  each $[\omega] \in \hat \Gamma_{2n}$
has $4n$ representatives in $\Gamma_{2n}$.
    Let 
$\hat \Gamma = \cup_{n \geq 0} \hat \Gamma_{2n}$ and
let $\hmusap$ be the measure on $\hat\Gamma$ that assigns
measure $\beta^{-n}$ to each $[\omega] \in \hat \Gamma_{2n}$.
Equivalently, we could assign measure $(4n)^{-1} \beta^{-n}$
to each $\omega \in \Gamma_{2n}$ and then define the
measure of $[\omega]$ to be the sum of the measures of
all its representatives.  In \cite{MS}, a self-avoiding
polygon is defined to be the set of equivalence
classes in $\Gamma_{2n}^*$, where two walks are equivalent
if a translation of one is equivalent to the other
under $\sim$. In particular, in their definition the number
of SAPs of length $2n$ is $(4n)^{-1} \; \#(\Gamma_{2n}^*)$.

\subsection{Critical exponents}

There are a number of ``critical exponents''
associated to SAWs and SAPs.  These are conjectured
to be {\em universal} quantities in that their
values will not change if one changes the lattice.
At the moment there is no proof that these exponents
exist.  We will introduce these
exponents in this subsection. 
We will not be rigorous in this section, e.g.,  we will
assume the existence of exponents.  However, this section
will be self-contained and   will
not assume any previously derived nonrigorous predictions
for the exponents.

\subsubsection{Mean-Square Displacement}

Let $\diam(\omega)$ denote the diameter of a SAW or SAP, i.e.,
the maximal value of $|\omega_j - \omega_k|$.  It
is believed that there exists an exponent $\nu$, called
the {\em mean-square displacement exponent}, such that
\[   [\#(\Lambda_n^*)]^{-1} \sum_{\omega \in
    \Lambda_n^*} \diam(\omega) \: \approx\:
    [\#(\Gamma_n^*)]^{-1} \sum_{\omega \in
       \Gamma_n^*} \diam(\omega) \:\approx\: n^{\nu} ,
  \;\;\;\; n \rightarrow \infty. \]
Recall that we are using $\approx$ to mean that the 
logarithms of both sides are asymptotic.  In fact, above and
in all the exponents below  it is believed that
the $\approx$ can be replaced with $\asymp$.
The exponent $\nu$ is often also described as in (\ref {nudef}), 
which should be roughly equivalent. 

\subsubsection{Number of walks}

The exponent $\gamma$ is defined by the relation
\[  \musaw(\Lambda_n^*) = \beta^{-n} \; \#(\Lambda_n^*)
               \approx n^{\gamma - 1} , \;\;\;\;
     n \rightarrow \infty. \]
Another way of viewing the exponent is to say
that the probability that two SAWs of length
$n$ starting at the origin do not intersect
(and hence can be joined to form a SAW of length
$2n$) is about $n^{1 - \gamma}$. Assuming the
existence of $\nu$ and $\gamma$, we can say that
the $\musaw$ measure of the set of $\omega
\in \Lambda^*$ of diameter between $R$ and
$2R$ looks like 
$$
 \sum_{n = R^{1/\nu}}^{(2R)^{1/ \nu}}
n^{\gamma -1} \asymp  R^{\gamma/\nu}
$$ as
$R \rightarrow \infty$.

\subsubsection{Number of half plane walks}  \label{halfsec}

Let $\Lambda_n^+$ denote the set of $\omega 
=[\omega_0,\ldots,\omega_n] \in
\Lambda_n^*$ such that
$\{\omega_1,\ldots,\omega_n\} \subset
  \{(x^1,x^2)\in \Z^2: x^2 > 0\}$, and
let $\Lambda^+ = \cup_{n \geq 0} \Lambda_n^+$.
The exponent $\rho$ is defined by
\[ \musaw(\Lambda_n^+)= \beta^{-n} \; \#(\Lambda_n^+) 
       \approx n^{\gamma - 1 - \rho}. \]
In other words, the probability that a SAW
chosen uniformly at random from $\Lambda_n^*$ lies in
the half plane looks like $n^{-\rho}$.     The
$\musaw$ measure of the set of $\omega 
\in \Lambda^+$ of diameter between $R$ and
$2R$ looks like
\begin{equation}  \label{1}
  \sum_{n= R^{1/\nu}}^{(2R)^{1/\nu}} n^{\gamma - 1
   - \rho} \asymp R^{(\gamma - \rho)/\nu} . 
\end{equation}

\subsubsection{Number of polygons}  \label{numpolysec}

The exponent $\alpha$ (which 
is denoted  $\alpha_{{\rm sing}}$ in \cite{MS})  is defined by
the relation
\[   \musap(\Gamma_{2n}^*) = \beta^{-2n}
           \; \#(\Gamma_{2n}^*) \approx
              n^{\alpha - 2} . \]
It is conjectured   that $\alpha -2 =   - 2\nu $.
We review the heuristic argument here.  A SAP
of $2n$ steps
consists of two SAWs of $n$ steps.  The
$\musaw \times \musaw$ measure on pairs of
$n$-step walks starting at the origin
is about $n^{2(\gamma - 1)}$.
Since the average diameter is about $n^{\nu}$,
we expect that the probability that the two
walks have the same endpoint is about $n^{-2\nu}$.
We also need the two walks not to intersect near
the origin (contributing a factor of $n^{1 - \gamma}$)
and not to intersect near their endpoints 
(contributing another factor of $n^{1 - \gamma}$).
Hence the {\em probability} that two walks chosen
uniformly from $\Lambda_n^* \times \Lambda_n^*$
can be joined to form a SAP should look like
$n^{-2\nu}n^{2 (1 - \gamma)}$ and hence the 
$\musap$ measure
of SAPs of $2n$ steps is about $n^{-2\nu}$.

If we assume this conjecture, then the $\musap$
measure of the set of SAPs starting at the origin
of diameter between
$R$ and $2R$  looks like
\begin{equation}  \label{3}
 \sum_{n=R^{1/\nu}}^{(2R)^{1/\nu}}
           n^{-2 \nu}  
     = R^{(1/\nu) - 2}. 
\end{equation}

  Let $\Gamma_{2n}^+$ be the set of
$\omega=[\omega_0,\ldots,\omega_{2n}] \in \Gamma_{2n}^*$
such that $\{\omega_1,\ldots,\omega_{2n-2}\} \subset
 \{(x^1,x^2)\in \Z^2: x^2 > 0\}$. Let $\Gamma^+ = \cup
\Gamma_{2n}^+$.
If $\omega
\in \Gamma_{2n}^*$, then either zero or two of the $4n$
representatives of $[\omega] \in \hat \Gamma_{2n}$ are in $\Gamma_{2n}^+$.
(In order for there to be more than zero, there must be
a unique horizontal bond with smallest $y$-coordinate.  If this bond
exists then there are two possible representatives, since
we can choose either vertex of this bond for the origin
and this choice determines the orientation.)
However, there is a one-to-one
map from $\hat \Gamma_{2n}$ into $\Gamma_{2n+2}^+$ given by considering
a horizontal edge of the SAP with least $y$-coordinate
(choosing arbitrarily if there is more than one),
replacing this edge by the three edges in the square face
of $\Z^2$ below it, and taking the representative in
$\Gamma_{2n+2}^+$ which has $\omega_{2n}=(1,0)$.
Since we believe that $\#\Gamma_{2n+2}^+/\#\Gamma_{2n}^+$
is bounded, we conclude that the probability that a
randomly uniformly chosen  $\omega
\in \Gamma_{2n}^*$ stays in the half plane is comparable to
$\# \hat \Gamma_{2n}/\#\Gamma_{2n}^*$, i.e., 
comparable to $n^{-1}$.  Hence,
\[   \musap(\Gamma_{2n}^+) \approx
              n^{\alpha - 3}  . \]
 The
$\musap$ measure of the set of $\omega 
\in \Gamma^+$ of diameter between $R$ and
$2R$ looks like
\begin{equation}  \label{2}
  \sum_{n= R^{1/\nu}}^{(2R)^{1/\nu}} n^{\alpha -3
   } \asymp R^{(\alpha -2)/\nu}  = R^{-2}. 
\end{equation}
We can see that SAPs are in some sense easier than
SAWs --- we have derived an exponent 2 without working too hard
and without
knowing exact values of any of the nontrivial exponents.

Note that $\hmusap(\hat \Gamma_{2n})
\approx n^{\alpha - 3} $. The $\hmusap$ measure of the set
of $[\omega]$ of diameter greater than $R$ contained in
a disk of radius $4R$, say, is $\approx 1$.
(There are $\approx R^2$ possible choices for the leftmost lowest point,
and fixing the leftmost lowest point, the measure is $\approx R^{-2}$
by the exponent derived above.)

\subsubsection{Nonintersecting walks} \label{nonintsec}

If $\omega = [\omega_0,\omega_1,\ldots,
\omega_n] \in \Lambda$ and $x \in \Z^2$, let
$\omega+x$ denote the translated walk
$[\omega_0 + x, \omega_1 + x, \ldots,\omega_n + x]$.
Let $e_1 = (1,0)$.
There are two natural exponents
associated with nonintersection of SAWs. 
    Suppose $n \leq m$
and $(\omega,\omega')$ is chosen from the uniform
distribution on 
  $\Lambda_n^* \times \Lambda_m^*$.  Note that  $\omega \cap (\omega' + e_1)
= \emptyset $ if and only if the ``reversal'' of $\omega$,
$[0,e_1]$ and $\omega'+ e_1$ 
can be concatenated  to
form a SAW. From this we can see that the
probability that $\omega \cap (\omega' + e_1)= \emptyset$
decays like $n^{1 -\gamma}$.  

If $(\omega,\omega')$
is chosen from the uniform
distribution on 
  $\Lambda_n^+ \times \Lambda_m^+$, then the probability
that $\omega \cap (\omega' + e_1) = \emptyset$ looks like
$n^{-1}$.  To see this,   note that
an argument similar to that in \S\ref{numpolysec} shows that a 
positive fraction of the SAWs should have a unique bond
with minimal second component.
In such a situation, the walk can be translated to give a
pair $(\omega,\omega')$ with $\omega \cap (\omega' + e_1) = \emptyset$.

\subsection{Walks with fixed endpoints}  \label{scalesec}

In this section we relate these exponents to 
those that describe 
the measure of sets of self-avoiding walks or
polygons with prescribed endpoints,
which will then 
become the ``scaling exponents'' in  
the scaling limit.    If $D$ is an open domain in
$\C  $, $x.y \in \Z^2$,
 and $\omega=[\omega_0,\ldots,\omega_n] \in \Lambda_n$,
  let $\Lambda(x,y;D,N)$
be the set of all $\omega =[\omega_0,\ldots,\omega_n]
\in \Lambda$ such that $\omega_0 = x, \omega_n =
y$ and   $[\omega_1,\ldots,\omega_{n-1}]
\subset ND$.  (When we write  $
[\omega_1,\ldots,\omega_{n-1}] \subset ND$. we
mean that  the bonds $[ \omega_j,
 \omega_{j+1}]$ are   in $ND$, $j=1,2,\ldots,n-2$.)
Let $\Gamma(x,y;D,N)$ be the set of all
$\omega = [\omega_0,\ldots,\omega_{2n}] \in
\Gamma$ such that $\omega_0 = \omega_{2n} = x$;
$y = \omega_k$ for some $k=1,2,\ldots,2n-1$; and
  $[ \omega_1,\cdots,\omega_{k-1}]  \cup
[\omega_{k+2},\ldots,\omega_{2n-2}] \subset
ND$. 
 (Note that we do not require
the bonds $ [\omega_{2n-2},\omega_{2n-1}]$ and
$ [\omega_{k+1},\omega_{k+2}]$ to be  
in $ND$; this allows the bonds $ [\omega_{2n-1},
\omega_{2n}]$ and $ [\omega_k,\omega_{k+1}]$ to be
entirely outside $ND$).
Each $\omega
\in   \Gamma(x,y;D,N)$  can be considered as a
pair of disjoint walks $\omega^1,\omega^2$
with $\omega^1  \in \Lambda(x,y;D,N),
\omega^2 \in \Lambda(x_1,y_1;D,N)$ and $|x - x_1| = 1,
|y - y_1|=1$.

\subsubsection{Boundary to boundary}

Let $D \subset \C$ be the open rectangle
$\{x + iy: -1 < x < 1, 0 < y < 1\}$ (here and 
in the next subsections, the particular choice of
domains $D$ will be somewhat arbitrary) and consider
$ \Lambda(0,Ni; D;N).$  This is
the lattice approximation, with lattice spacing
$N^{-1}$, of self-avoiding paths from boundary
point $0$ to boundary point $i$ in
$D$.
The {\em SAW boundary scaling exponent} $a$ is defined by
the relation
\[  \musaw[\Lambda(0,Ni;D,N)] 
     \approx N^{-2a} , \;\;\;\; N \rightarrow \infty. \]
Let us give a heuristic guess as to what this
exponent should be in terms of the critical
exponents $\nu,\gamma,\rho$.  By (\ref{1}),
the $\musaw$ measure
of the set of SAWs starting
at the origin of diameter
of order $N$ (e.g., between $N$
and $2N$) whose second
component stays positive
is about $N^{(\gamma - \rho)/\nu}$.  The probability
that the final point of such a walk is
$Ni$ should look like $N^{-2}$ (since there
are $O(N^2)$ points distance $O(N)$ from
the origin).  We also need the  the walk to
stay below the line $Ni +\Z $;
 the probability
that this happens should look like
$N^{-\rho/\nu}$ (remember the walk has about
$N^{1/\nu}$ steps).  Using this, we get the
conjecture
\begin{equation}  \label{oct27.1}
         a = 1 + \frac{2 \rho - \gamma}{2\nu}. 
\end{equation}

The {\em SAP boundary scaling exponent} $a'$ is defind by
the relation
\[  \musap[\Gamma(0,Ni;D,N)] 
     \approx N^{-2a'} . \]
By (\ref{2}), the $\musap$ measure of the
set of SAPs starting at the origin of diameter of
order $N$ whose second component stays positive
is about $N^{-2}$.  A positive fraction of these walks
should have a unique bond 
  with maximal
second    component. The
probability that one of the vertices of this bond
is $Ni$ is
of order $N^{-2}$.  From this we get the
conjecture  
\[  \musap[\Gamma(0,Ni;D,N)] 
     \approx  N^{-4} ,\]
i.e., $a'=2$.
 
\subsubsection{Boundary to interior}  \label{boundintsec}

Let $D$ be the open square 
$\{x + iy: -1 < x < 1, -1 < y < 1\}$ and consider
$ \Lambda(0,Ni; D;N).$  This is
the lattice approximation, with lattice spacing
$N^{-1}$, of self-avoiding paths from interior
point $0$ to boundary point $i$ in
$D$. 
 Define the
{\em SAW interior scaling exponent} $b$ by 
\[  \musaw[\Lambda(0,Ni;D,N)] 
     \approx  N^{-(a+b)} . \]  
The heuristics for $b$ are done similarly
as for $a$.  In this case 
the $\musaw$ measure
of the set of SAWs starting
at $0$ of diameter
of order $N$  
is   $N^{\gamma  /\nu}$.  Note that
we are not requiring the walk to stay
in the upper half plane so this factor
 differs from that in the previous case.
However, the probability of having
$Ni$ as the terminal point and having
the walk stay below the line $y = N$ contribute
factors of $N^{-2}$ and $N^{-\rho/\nu}$ as before.
Hence, we get the conjecture 
\[  a+ b = 2 + \frac{ \rho - \gamma}{\nu}, \]
or
\begin{equation}  \label{oct28.2}
  b = 1 - \frac{\gamma}{2\nu} . 
\end{equation}  

The {\em SAP interior scaling exponent} $b'$ is defind by
the relation
\[  \musap[\Gamma(0,Ni;D,N)] 
     \approx  N^{-(a'+b') }  = N^{-(2 + b')}. \]
By (\ref{3}), 
the $\musap$ measure of the set of SAPs
of diameter of order $N$ (i.e., between
$N$ and $2N$, say)
looks like $N^{(1/\nu) - 2}$.
Again, the probability that $Ni$ is a point
with maximal second component is of order $N^{-2}$.
This gives the conjecture
\[      b' = 2 - \frac{1}{\nu}. \]

\subsubsection{Interior to interior}

Let $D \subset \C$ be the open rectangle
$\{x + iy: -1 < x < 1, -1 < y < 2\}$ and consider
$ \Lambda(0,Ni; D;N).$  This is
the lattice approximation, with lattice spacing
$N^{-1}$, of self-avoiding paths from interior
point $0$ to interior point $i$ in
$D$.
  The probability
that a randomly chosen SAW starting at $0$
of diameter between $N$ and $2N$ stays in $D$
and ends at $Ni$ should be of order $N^{-2}$.
  Hence, we expect that 
\[  \musaw[\Lambda(0,Ni;D,N)] 
     \approx N^{(\gamma/\nu) - 2} 
   = N^{-2b} . \]
Similarly,
\[ \musap[\Gamma(0,Ni;D,N)]
   \approx N^{-2b'} = N^{(2/\nu) - 4} . \]

\subsection{Scaling limit}

\subsubsection{Convergence of measures on curves}

There are a number of   ways to define convergence
of measures on curves in $\R^2 = \C$.   A standard way is to define
a metric  on the space of curves and then to use
weak convergence of measures on the generated
Borel $\sigma$-algebra.   The metrics we define are
actually pseudometrics, i.e., metrics on equivalence classes
of curves where two curves are equivalent if their distance
is zero.

Suppose $\gamma^1:[0,T^1] \rightarrow \C$,
$\gamma^2:[0,T^2] \rightarrow \C$ are continuous curves.
Define the Hausdorff metric between the curves by
\[  d_H(\gamma^1,\gamma^2) = 
    \max\bigl\{  \; \max\{ \dist(z,\gamma^2[0,T^2]);
             z \in \gamma^1[0,T^1]\} ,  \hspace{1in} \]
\[  \hspace{1.5in}
     \max\{\dist(z,\gamma^1[0,T^1]); z \in \gamma^2[0,T^1] \} \; \bigr\}
   . \]
The Hausdorff metric depends only
on  
the sets $\gamma^1[0,T^1]$ and $\gamma^2[0,T^2]$.

A finer topology is obtained by using the metric $d$ where
$d(\gamma^1,\gamma^2) < \epsilon$ if there is an increasing
homeomorphism $\phi:[0,T_1] \rightarrow [0,T_2]$ such that
\[    \bigl|\gamma^2(\phi(t)) - \gamma^1(t) \bigr| < \epsilon,
  \;\;\;\; 0 \leq t \leq T_1 . \]
Note that $d(\gamma^1,\gamma^2) \geq d_H(\gamma^1,
\gamma^2)$ and $d(\gamma^1,\gamma^2) = 0$ if and
only if $\gamma^2$ can be obtained from $\gamma^1$ by
a monotone nondecreasing reparametrization.  This metric considers
not just the set of points visited by a curve but also the
order in which they are visited; however, it does not
consider the ``speed'' at which the curves are
traversed.

Call $\gamma^1, \gamma^2$ {\em closed} curves
if $\gamma^1(T_1) = \gamma^1(0)$ and $\gamma^2(T^2)
= \gamma^2(0)$.   Let $[0,T_1]^*$ denote the topological
circle obtained by identifying $0$ and $T_1$.  The
metric $d_C$ on closed curves is defined by saying
that 
$d_C(\gamma^1,\gamma^2) < \epsilon$ if there is
a homeomorphism $\phi:[0,T_1]^* \rightarrow
[0,T_2]^*$ such that 
\[    |\gamma^2(\phi(t)) - \gamma^1(t) | < \epsilon,
  \;\;\;\; 0 \leq t \leq T_1 . \]
(In our definition we allow the $\phi$ to have
either orientation; we could   define a metric
$d_O$ similarly
by restricting to positively oriented homeomorphisms $\phi$.)

We will make conjectures about convergence of measures.
We will not be explicit about the form of the conjectures;
implicitly we will mean convergence using the metric
$d$ or $d_C$.   However, proofs of convergence using
the weaker metric $d_H$ would be interesting and many
of the conclusions probably need  only $d_H$ convergence.
 In any case, we will  be interested in curves only
{\em up to reparametrization}; for this reason, we will
not be explicit about the parametrization of scaled 
SAWs and SAPs.  We may always assume that we
have chosen a parametrization of $\gamma^1$ with
$T_1 = 1$.

Let us stress now that when  we  later 
conjecture that $SLE_{8/3}$ and the
scaling limit of SAW have
the same distribution, this should be considered only
as a statement about curves modulo parametrization:
The ``natural'' parametrization that one would give
to the scaling limit of SAWs (which is 
 related to the dimension of the curve) 
is {\em not} the same
as the parametrization of $SLE_{8/3}$ defined in the 
previous section.

\subsubsection{Scaling limits for SAWs}  \label{sawscalesec}

Let $D$ be a connected open domain; for ease
 we will assume that $D$ is bounded.  Let
$z,w$ be distinct points in $D$.   
  For each positive integer $N$,
let $z_N,w_N$ be lattice points closest
to $Nz,Nw$   (We need to make some convention if 
there is more than one    $z_N,w_N$;
we can make any choice.)
Assume that the following stronger
version of the scaling law holds: there is
a $C(z,w;D) \in (0,\infty)$ such that
\[  \lim_{N \rightarrow \infty}
     N^{2b} \; \musaw[\Lambda(z_N,w_N;D,N)]
              = C(z,w;D) . \]
To each $\omega \in \Lambda(z_N,w_N;D,N)$, we
consider $N^{-1} \omega$ as a continuous curve
connecting $N^{-1} z_N$ to $N^{-1} w_N$; as we
have mentioned before, we may use any parametrization
of this curve.
 We conjecture that as $N \rightarrow
\infty$, the measure $N^{2b}\; \musaw$,
restricted to $\Lambda(z_N,w_N;D,N)$,
considered as a measure on curves in $D$ by the relation
$\omega \mapsto  N^{-1} \omega $, converges weakly to
a measure $\msaw(z,w;D)$.  By construction, the
total mass of this measure is $C(z,w;D)$.

 Suppose that $\p D$ is smooth,  
$z \in D$, and let $z_N$ be as above.  Let
$\bar \Lambda(z_N;D,N)$ be the set of
SAWs $ [\omega_0,\ldots,\omega_n]$ with
$\omega_0 = z_N, [\omega_0,\ldots,\omega_{n-1}]
\subset ND$, $[\omega_{n-1},\omega_n] \cap
\p D \neq \emptyset$. 
We assume there is a constant $C(z;D) \in
(0,\infty)$ such that
\[  \lim_{N \rightarrow \infty} N^{a+b}\;
      N^{-1} \; \musaw[\bar \Lambda(z_N;D,N)]
    = C(z;D) . \]
The $N^{-1}$ appears because there are $O(N)$ 
``boundary points'' at grid scale $1/N$.
We conjecture that as $N \rightarrow
\infty$, the measure $N^{a+b}\;
N^{-1} \;  \musaw$, restricted
to $\bar \Lambda(z_N;D,N)$,
considered as a measure on curves in $\overline D$ by the relation
$\omega \mapsto N^{-1} \omega$, converges weakly to
a measure $\msaw(z;D)$.  By construction, the
total mass of this measure is $C(z;D)$.  We
also assume that this measure can be written as
\[   \msaw(z;D) = \int_{\p D} \msaw(z,w;D) \;
    d|w| , \]
where $\msaw(z,w;D)$ is a measure on paths connecting
the interior point $z$ to the boundary point $w$.  We
let $C(z,w;D)$ be its total measure so that
\[   C(z;D) = \int_{\p D} C(z,w;D) \; d|w| . \] 
Of course, this defines $C(z,w;D)$ only for almost every
$w\in\partial D$.  However, we assume that $C(z,w;D)$ may
be chosen to be continuous in $w$.

Similarly, we can consider measures connecting
two distinct points $z,w \in \p D$. 
  Let $A_1,A_2$ be nontrivial, disjoint,
closed  subarcs of $\p D$.
Let $\musaw(A_1,A_2;D,N)$ be the set of all
$\omega = [\omega_0,\omega_1,\ldots,\omega_n] \in \Lambda$ such that
\begin {itemize}
\item
$[\omega_1,\ldots,\omega_{n-1}]
\subset ND, $
\item 
 $[\omega_0,\omega_1] \cap N A_1 \neq \emptyset$ and $
[\omega_{n-1},\omega_n] \cap N A_2\neq \emptyset$. 
\end {itemize}
Assume that  there is
a $C(A_1,A_2;D) \in (0,\infty)$ such that
\[  \lim_{N \rightarrow \infty}
     N^{2a} \; N^{-2} \; \musaw[\Lambda(A_1,A_2;D,N)]
              = C(A_1,A_2;D) . \]
We conjecture that as $N \rightarrow
\infty$, the restriction of the 
measure $N^{2a-2}\musaw$,
to $\Lambda(A_1,A_2;D,N)$,
   converges weakly to
a measure $\msaw(A_1,A_2;D)$.   By construction, the
total mass of this measure is $C(A_1,A_2;D)$.  
We also assume that this measure can be written as
\[  \msaw(A_1,A_2;D) = \int_{A_1} \int_{A_2}
    \msaw(z,w;D) \; d|w|\; d|z|  , \]
where $\msaw(z,w;D)$ is a measure on curves connecting
boundary points $z$ and $w$ with total mass $C(z,w;D)$
satisfying
\[     C(A_1,A_2;D) = \int_{A_1} \int_{A_2} C(z,w;D)
    \; d|w|\; d|z|. \]

  Note that the meaning of $\msaw(z,w;D)$
varies whether both points are interior points;
both points are boundary points; or one is an
interior point and one is a boundary point.
Hopefully, this should not cause confusion.

   
\subsubsection{Conformal covariance}
The major
 assumption we will make about the limiting
measure (other than its existence as a scaling limit
of discrete approximations)
is {\em conformal covariance}.  Suppose, for ease, that
$D$ is a bounded domain with smooth boundary $\p D$
and $f: D \rightarrow D'$ is a conformal transformation
such that $f$ is $C^1$
up to the boundary.
   The conformal
covariance can be phrased as follows.

\begin{itemize}
\item If $z,w \in D$, then
\[    f \circ \msaw(z,w;D) = |f'(z)|^b \;
         |f'(w)|^b \; \msaw(f(z),f(w);D') . \]
\item If $z \in D, w \in \p D$,
 \[    f \circ \msaw(z,w;D) = |f'(z)|^b \;
         |f'(w)|^a \; \msaw(f(z),f(w);D') . \]
\item If $z , w \in \p D$,
 \[    f \circ \msaw(z,w;D) = |f'(z)|^a \;
         |f'(w)|^a \; \msaw(f(z),f(w);D') . \]
\end{itemize}
In other words, the image of the measure under
a conformal transformation is the same as the measure
in the new domain  multiplied by a scalar correction
factor.   
The scaling laws tells us that these relations 
hold for the maps $z \mapsto rz$ if $r > 0$.   
The exponents $a,b$ are sometimes called {\em conformal
weights}.

\subsubsection{Scaling limit for SAPs}  \label{sapscalesec}

We can make similar assumptions about the limit of SAPs.
Since the definitions are virtually the same as for
SAWs, we will not repeat them.  These gives us measures
$\msap(z,w;D)$ on closed curves going through both
$z$ and $w$ lying in $D \cup \{z,w\}$.  
We also assume that these curves have
the property that $z,w$ are on their boundary, i.e., the
curves do not disconnect $z$ from $w$.  

We also assume that these curves are conformally covariant.
Recall that $a' = 2, b' = 2 - (1/\nu)$. If $f,D,D'$ are
as in the previous subsection, we assume
\begin{itemize}
\item If $z,w \in D$, then
\[    f \circ \msap(z,w;D) = |f'(z)|^{2 - (1/\nu)} \;
         |f'(w)|^{2 - (1/\nu)} \; \msap(f(z),f(w);D') . \]
\item If $z \in D, w \in \p D$,
 \[    f \circ \msap(z,w;D) = |f'(z)|^{2 - (1/\nu)} \;
         |f'(w)|^2 \; \msap(f(z),f(w);D') . \]
\item If $z , w \in \p D$,
 \[    f \circ \msap(z,w;D) = |f'(z)|^2 \;
         |f'(w)|^2 \; \msap(f(z),f(w);D') . \]
\end{itemize}

\subsubsection{Restriction property}  \label{restsec}
\label {secheuristic}

If $z,w \in D$, write $\musaw(z,w;D,N)$ for 
$N^{2b} \; \musaw$ restricted to 
$\Lambda(z_N,w_N;D,N)$.  By
definition, the measures $\musaw(z,w;D,N)$
satisfy the {\em restriction property}:
if $D' \subset D$, then $\musaw(z,w;D',N)$ 
is $\musaw(z,w;D,N)$ restricted to walks
that stay in $N D'$ (in the appropriate sense).
Hence,  the limit measure $\msaw(z,w,D)$,
assuming it exists,
  must satisfy this property.   Similarly, this
property should hold in the limit if $z \in D, w \in
\p D$ or $z,w \in \p D$.
In other words, all the families
$\msaw$ and $\msap$ should be satisfy restriction.

\medskip

It now follows that
the limit $\msaw$ is supported on curves $\gamma:
[0,1] \rightarrow \C$ with $\gamma[0,1] \subset D\cup\{z,w\}$.  
(The points $z$ and $w$ may be in $D$ or in $\partial D$.)
This is because when $D$ is smooth we may find
a smooth domain $D'\subset D$ such that $\partial D' \subset D\cup\{z,w\}$
and the triple $(D',z,w)$ is close, in the appropriate sense,
to the triple $(D,z,w)$.  We expect the total mass $C(z,w;D)$ of $\msaw(z,w;D)$
to be continuous in the triple $(D,z,w)$.
Consequently, most of the mass of $\msaw(z,w;D)$ is in the support
of $\msaw(z,w;D')$, and therefore is on paths in
$D\cup\{z,w\}$.

This also helps to understand why the limit 
measure should lie on simple curves.
If we condition on an initial segment of the SAW,
for the rest of the SAW, that initial segment acts like a boundary.
Since the curve avoids the boundary, it should also avoid its past;
i.e., it should not have any
self-intersections.

\subsubsection{Infinite SAW and SAP}

The infinite SAW and SAP are processes that are supposed
to look like the ``beginning'' of a long SAW or SAP.  
Their definition depends on the existence of a limit
that has not been proved to exist in two dimensions
except in the case of the half-plane infinite SAW.

If $\omega
=[0,\omega_1,\ldots,\omega_k] 
\in \Lambda^*, \omega^1 = [0,\omega_1^1,\ldots
\omega_n^1] \in \Lambda^*$, define the concatenation
$\omega \oplus \omega^1 = [0,\omega_1,
\ldots,\omega_k,\omega_k + \omega_1^1,\ldots,$ $
\omega_k + \omega^1_n]$.  It is possible that
the concatenation is not self-avoiding.
 Let $\Lambda_n^*(\omega)$ be the set
of $\omega^1 \in \Lambda_n^*$ such that $\omega \oplus
\omega^1 \in \Lambda^*$, and if $\omega \in \Lambda^+$,
let $\Lambda^+_n(\omega)$ be the set of
 $\omega^1 \in \Lambda_n^*$ such that $\omega \oplus
\omega^1 \in \Lambda^+$.    For $\omega \in \Lambda_k^*$,
$\omega' \in \Lambda_k^+$, let 
\[  Q(\omega) = \lim_{n \rightarrow \infty}
    \frac{\#[\Lambda_n^*(\omega)]}{\#[\Lambda_{n+k}^*]} ,
  \;\;\;\;  Q^+(\omega') = \lim_{n \rightarrow \infty}
    \frac{\#[\Lambda_n^+(\omega')]}{\#[\Lambda_{n+k}^+]} , \]
assuming the limit exists.   Roughly speaking, $Q(\omega)$
is the fraction of very long SAWs whose beginning
is $\omega$, and similarly for $Q^+(\omega')$. The limit
$Q(\omega)$ is not known to exist, although it
is believed to.   In the appendix 
we adapt an argument of Madras-Slade \cite {MS} to 
  prove that the limit $Q^+(\omega')$ exists.

By convention, let us define $Q(\omega) = 0$ if
$\omega \not\in \Lambda$ and $Q^+(\omega') = 0$
if $\omega' \not\in \Lambda^+$.  Then, assuming the
existence of the limits, for each $\omega$ and
each positive integer $n$,
\[     Q(\omega) = \sum_{\omega' \in \Lambda_n^*}
                 Q(\omega \oplus \omega') , \;\;\;\;
     Q^+(\omega) = \sum_{\omega' \in \Lambda_n^*}
                 Q^+(\omega \oplus \omega'). \]
The {\em infinite self-avoiding walk} starting at
the origin is the stochastic process $X_j$ 
such that for all $\omega = [0,\omega_1,
\ldots,\omega_k] \in \Lambda^*$,
\[ \Prob\{X_0 = 0, X_1 = \omega_1, \ldots,
   X_k = \omega_k\} = Q(\omega) . \]
Of course, the existence of this process depends
on the existence of $Q$.
Similarly the {\em half plane infinite self-avoiding
walk} starting at the origin is the stochastic
process $X_j'$ with
\[ \Prob\{X_0' = 0, X_1' = \omega_1, \ldots,
   X_k' = \omega_k\} = Q^+(\omega) . \]
  Assuming the existence
of the limit,  $Q$ 
is also equal to the limit (for $\omega
\in \Lambda_k^*$)
\[   Q(\omega) = \lim_{r \rightarrow (1/\beta)-}
     \frac{\sum_{n=0}^\infty  r^{n+k} \; \#[\Lambda_{n+k}^*(\omega)]}
    { \sum_{n=0}^\infty  r^{n+k}\; \#[\Lambda_n^*]} , \]
   We cannot put $r = 1/\beta$
directly into the sums, because the sums  
 diverge \cite[Corollary 3.1.8]{MS}.   It is possible that the
limit above could be easier to establish than the limits in the 
definition of $Q$, and that would give a reasonable alternative definition
for the infinite SAW.

An infinite SAP going through the
origin is a pair of infinite SAWs
that do not intersect.  Let $\Lambda_k^{*,2}$ be the
set of ordered pairs $(\omega^1,\omega^2) \in
\Lambda_k^* \times \Lambda_k^*$ with $\omega^1 \cap
\omega^2 = \{0\}$.
Similarly, let $\Lambda_k^{+,2}$ be the
set of ordered pairs $(\omega^1,\omega^2) \in
\Lambda_k^+ \times \Lambda_k^+$ with $\omega^1 \cap
(e_1+\omega^2) =\emptyset$.  If $(\omega^1,\omega^2) \in
\Lambda_k^{*,2}$, let $\Lambda_n^{*,2}(\omega^1,\omega^2)$ 
be the set of $(\omega^3,\omega^4) \in
\Lambda_n^* \times \Lambda_n^*$ with
$(\omega^1 \oplus \omega^3,\omega^2 \oplus \omega^4)
\in \Lambda_{k+n}^{*,2}$.  Define $\Lambda_n^{+,2}
(\omega^1,\omega^2)$ similarly (for $(\omega^1,
\omega^2) \in \Lambda_k^{+,2}$), and let
\[  Q(\omega^1,\omega^2) = \lim_{n \rightarrow \infty}
    \frac{\#[\Lambda_n^{*,2}(\omega^1,
  \omega^2)]}{\#[\Lambda_{n+k}^{*,2}]} ,
  \;\;\;\;  Q^+(\omega^1,\omega^2) = \lim_{n \rightarrow \infty}
    \frac{\#[\Lambda_n^{+,2}(\omega^1,\omega^2
   )]}{\#[\Lambda_{n+k}^{+,2}]} . \]
The {\em infinite self-avoiding polygon} is a
pair of stochastic processes $(X^1_j,X^2_j)$ such
that for each $(\omega^1,\omega^2) \in
\Lambda_{k}^{*,2}$,
\[  \Prob\{X^1_0= \omega_0^1,\ldots,X^1_k= \omega_k^1;\;
   X^2_0 = \omega_0^2,\ldots, X^2_k = \omega_k^2\}
    = Q(\omega^1,\omega^2). \]
The {\em  half plane infinite self-avoiding polyon}
is defined similarly using $Q^+$.

\subsubsection{Scaling limit of infinite SAW and SAP}

The scaling limits of SAWs and SAPs have been written in
terms of finite measures that are not probability measures.
Sometimes, it is more convenient to consider the associated
probability measures 
 \[ \msaw^\#(z,w;D) = \frac{\msaw(z,w;D)}{  |\msaw(z,w;D)|},
\;\;\;\;
 \msap^\#(z,w;D) = \frac{\msap(z,w;D)}{  |\msap(z,w;D)|}. \]
The conformal covariance  assumption implies
that the measure families $\msaw^\#(z,w;D)$ and $\msap^\#(z,w;D)$
are conformally {\em invariant}.
This 
conjecture can be extended to cases where the boundary
of $D$ is not smooth. 

For the boundary to boundary case,
assuming conformal invariance, it suffices to determine
the distribution of $\msaw^\#(0,\infty;\H)$ and
$\musap^\#(0,\infty;\H)$.  We conjecture
that these measures are the scaling limits of
the half plane infinite SAW and the half plane
infinite SAP, respectively.  Similarly, in the
interior to interior case, we conjecture that
$\msaw^\#(0,\infty;\C)$ and $\msap^\#(0,\infty;
\C)$ are the scaling limits of the infinite SAW and the infinite SAP,
respectively.

\subsubsection{Boundary polygons}  \label{boundpolysec}

Let $\bar \Gamma(D;N)$ be the set of SAPs
$\omega = [\omega_0,\omega_1,\ldots,\omega_{2n}]$
with $[\omega_1,\ldots,\omega_{2n-2} ] \subset
ND$ and $[\omega_0,\omega_1] \cap N\p  D \neq
\emptyset$.
   For every $\epsilon > 0$ consider
the measure $\musap$ restricted to
 $\bar \Gamma(D;N;\epsilon),$ defined to be the set of
$\omega \in \bar \Gamma(D;N)$ such that
$\diam(\omega) \geq \epsilon N$.  By (\ref{2}),
we expect that
\[     \musap[\bar \Gamma(D;N;\epsilon)] 
    \approx c(D;\epsilon) N^{-1} . \]
(There is an extra factor of 
 $N $ since there are 
$O(N)$ boundary points).  We conjecture that
the measures  
$N \;\musap$ restricted
to $\bar \Gamma(D;N;\epsilon)  $
converge weakly (after rescaling the paths by a factor of $N^{-1}$) to a measure
$\mbarsap(D;\epsilon)$.  This is a measure on
curves of diameter at least $\epsilon$. We define
$\mbarsap(D)$ to be the measure on all curves
of positive diameter obtained as the limit
as $\epsilon \rightarrow 0+$ of $\mbarsap(D;\epsilon)$.
This is an infinite measure.  (We have defined
$\mbarsap(D)$ by means of the $\mbarsap(D;\epsilon)$ measures
in order to obtain a measure that gives zero mass to curves
of diameter $0$.  If we had taken the limit
of the measures $N \; \musap$ restricted
to all of $\Gamma(D;N) $, the limit
would give infinite measure to trivial curves.)

We conjecture that one can write
\[   \mbarsap(D) = \int_{\p D} \mbarsap(w;D) \; 
    d|w| , \]
where $\msap(w,D)$ is the appropriate measure
on boundary loops that hit $\p D$ at $w$.
This would define $\mbarsap(w;D)$ for almost all $w$ in $\partial D$.
We believe that one can choose a definition of $\mbarsap (w;D)$ for
all $w\in\partial D$ so that:
\begin{itemize}
\item $\mbarsap(w;D)$ is an infinite measure, but
the measure of curves of diameter at least $\epsilon$
is finite;
\item Conformal covariance: if $f:D \rightarrow
D'$ is a conformal transformation, then
\[   f \circ \mbarsap(w;D) = |f'(w)|^2
    \; \mbarsap(f(w);D') ;\]
(The scaling exponent $2$ comes from (\ref{2}).)
\item Restriction property: if $D' \subset D$,
and $w \in \p D' \cap \p D$,
then $\mbarsap(w;D')$ is $\mbarsap(w;D)$
restricted to curves that stay in $D'$.
\end{itemize}

We can also conjecture: $\mbarsap(w;D)$ is
a measure on simple loops, i.e.,
homeomorphic images of the unit circle, which touch
$\p D$ only at $w$.

\subsubsection{Soup of SAPs}  \label{sapsoupsec}

Recall  the $\hmusap$ measure on $\hat \Gamma$
from \S\ref{sawdefsec} and \S\ref{numpolysec}.  The scaling exponent for
that measure
is zero.  If $D$ is a domain and $N$ is a positive
integer, let $\hat \Gamma(D,N)$ be the set of 
$[\omega] \in \Gamma$ such that $N^{-1} \omega$ lies
entirely in $D$; note that this condition is independent
of the choice of representative $\omega$.  We conjecture
that the measure $\hmusap$ on $\hat \Gamma(D,N)$,
 on curves  converges weakly (when the curves are rescaled
by a factor of $N^{-1}$ to an (infinite)
measure $\mhsap(D)$.  This is an infinite measure
because it gives large measure to curves
of small diameter.  
However,
if $D$ is bounded and
$z,w \in D$ are distinct, then the $\mhsap(D)$ measure
of curves that ``surround'' both $z$ and $w$ should
equal some  $\hat c(D;z,w)< \infty$.  Let $\mhsap(D;z,w)$
denote $\mhsap(D)$ restricted to curves that surround
$z$ and $w$.

We make the assumption that $\mhsap(D)$ gives
no measure to 
trivial curves. (If necessary, we define the limit
in a way so that such curves are discarded.)
 
We conjecture that this measure satisfies
conformal {\em invariance}: if $f:D \rightarrow D'$ is
a conformal transformation, then $f \circ \mhsap(D)
= \mhsap(D')$.  Moreover, if $z,w \in D$, then
$f \circ \mhsap(D;z,w) = \mhsap(D';f(z),f(w))$.
Note that we are assuming
that there  is no scaling factor --- this is
a consequence of the fact that the scaling exponent
is zero.

The construction of this measure (assuming
existence of the limit) implies that $\mhsap$
satisfies the restriction property: if $D'
\subset D$, then $\mhsap(D)$ restricted to sets
staying in $D'$ is $\mhsap(D')$.

Note that both $\mbarsap$ and $\mhsap$ can be
considered as measures on hulls by considering the
hull bounded by the curves.

\section{Relating $SLE_\kappa$ to SAW and SAP}  \label{slesawsec}

In this section we discuss conjectures relating
SAWs and SAPs to $SLE_{8/3}.$ For some of the results,
we could have made theorems of the type, ``if the SAW
scaling limit exists and is conformally covariant then
the scaling limit is $SLE_{8/3}$.''  Since the
purpose of this paper is primarily to present conjectures,
we have not checked 
 the exact
hypotheses under which such conditional theorems could be attained.

\subsection{Boundary to boundary}

Consider the conjectured scaling limit $\msaw(z,w;D)$
as described in \S\ref{sawsec} where $D$ ranges
over smooth simply connected domains $D$ and $z,w \in
\p D$.  We can summarize the conjectures about
$\msaw(z,w;D)$ in the terminology of \S\ref{slesec2}.
Let   \[ \msaw^\#(z,w;D) = \frac{\msaw(z,w;D)}{  |\msaw(z,w;D)|}\]
be the probability measure associated to $\msaw(z,w;D)$.
Then  $\msaw^\#(z,w;D)$  is conjectured to be a $\rest(a)$
family.  The only $\rest(a)$ family supported on simple
curves is $a = 5/8$.  Recall that in this case
$\msaw^\#(z,w;D)$ is given by chordal
$SLE_{8/3}$ connecting $z$ and $w$ in $D$.

Another way of justifying the conjecture is to recall
that $\msaw^\#(0,\infty;D)$ is conjected to be the
scaling limit of the half plane infinite SAW.  This process
should have the conformal invariance properties discussed
in \S\ref{chordslesec}.  However, it should also
satisfy the restriction property.  This leaves $\kappa = 8/3$
as the only candidate.

\begin{conclusion}
The scaling limit of the half plane infinite SAW is
chordal $SLE_{8/3}$.  Moreover, the family of
measures
$\msaw^\#(z,w;D)$, where $D$ ranges over simply
connected domains and $z,w \in \partial D$,
is the unique $\rest(5/8)$
family.  In particular, the SAW boundary scaling
exponent is $a=5/8$.
\end{conclusion}

Recall from \S\ref{chordslesec}, that the Hausdorff
dimension of $SLE_{8/3}$ paths is $4/3$.  But the
dimension of the paths of the scaling limit of
SAWs should be $1/\nu$. Hence we get the following.

\begin{conclusion}
The mean-square displacement exponent
for SAWs and SAPs is  $\nu = 3/4$.
\end{conclusion}

Similarly, the family of measures
$\musap(z,w;D)$ are conjectured to give the 
$\rest(2)$ family.  This family is given by the hull
of the union of two independent Brownian excursions going
from $z$ to $w$ in $D$.

\begin{conclusion}
The scaling limit of the half plane infinite SAP is
the outer boundary of the union of two Brownian
excursions going from $0$ to $\infty$ in $\H$.
Moreover, the family of measures  
$\msap^\#(z,w;D)$,  where $D$ ranges over simply
connected domains and $z,w \in \partial D$,
is the unique $\rest(2)$
family.
\end{conclusion}

\subsection{Boundary to interior}

If we consider $\msaw^\#(z,w;D)$ where $z$ is an interior
point and $w$ is a boundary point, then a similar argument
allows us to conjecture that this measure comes from
a {\em radial} $SLE_\kappa$.  Since the dimension of the
paths should be the same as for the boundary to boundary
case, the only possible $\kappa$ is $8/3$.  

\begin{conclusion}
The family of
measures
$\msaw^\#(z,w;D)$, where $D$ ranges over simply
connected domains and $z \in D ,w \in \partial D$,
is  given by radial $SLE_{8/3}$ from $w$ to
$z$ in $D$.   
\end{conclusion}

  The following
prediction uses the rigorous scaling rule for
radial $SLE_{8/3}$ (\ref {548}), together with
(\ref{oct27.1}) and (\ref{oct28.2}).

\begin{conclusion}
The SAW interior scaling exponent is $b = 5/48$.
Also, if $\gamma,\rho$ are as defined in   
\S\ref{sawsec},
\[  \gamma = \frac{43}{32}, \;\;\;\;\;
  \rho = \frac{25}{64} . \]
\end{conclusion}

The measures $\msap^\#(z,w;D)$, where $D$ ranges over simply
connected domains and $z \in D ,w \in \partial D$,
are conjectured to satisfy the restriction property
and to be conformally covariant with boundary
scaling exponent $a' = 2$ and interior scaling
exponent $b' = 2 - (1/\nu)$.  This can be considered as
a measure on hulls by considering the hull surrounded
by the curves.  The prediction $\nu =
4/3$, gives the prediction $b'=2/3$.  It is not 
difficult to show that there is only one measure on
hulls that satisfies the restriction property
and is conformally covariant with these exponents.
In section \ref {S.radial}, we described this family
of measures in terms of non-disconnecting Brownian 
motions.

\begin{conclusion}
The family of
measures
$\msap^\#(z,w;D)$, where $D$ ranges over simply
connected domains and $z \in D ,w \in \partial D$,
is  given  by the outer boundary of a pair of
Brownian paths starting at $w$ conditioned to 
leave $D$ at $z$ and conditioned not to disconnect
$w$ from $\p D$.
\end{conclusion}

\subsection{Interior to Interior}

Consider the measures $\msaw^\#(0,\infty; \C)$
and $\msap^\#(0,\infty, \C)$.  These are the
conjectured scaling limits of the infinite
SAW and the infinite SAP, respectively.  

\begin{conclusion}
The scaling limit of the infinite SAW is the whole
plane $SLE_{8/3}.$
\end{conclusion}

\begin{conclusion}
The scaling limit of the infinite SAP is the measure
given by the outer boundary of a pair of non-disconnecting
Brownian motions starting at the origin.
\end{conclusion}

\subsection{Boundary bubbles} 

We end by conjecturing that the scaling of limit
of boundary polygons considered as a measure
on hulls is a scalar multiple of the measure on hulls
obtained from  boundary loops of Brownian motion
as described in \S\ref{conseqsec}, 
and the scaling limit of the soup of SAPs 
as a measure on hulls is a multiple of
the unrooted Brownian loop measure considered
as a measure on hulls.

\subsection {Other exponents}

Similar arguments based on the derivation of
the corresponding critical exponents for $SLE_{8/3}$
also yield the value of the exponents describing the 
non-intersection probabilities of self-avoiding walks 
(chordal or radial) that had been conjectured
by Duplantier and Saleur \cite {DS}, or those results 
concerning the exponents associated to the ``multi-fractal
spectrum'' of SAW (see \cite {Dmulti}); see also \cite {LW}.

\section*{Appendix: The infinite half-space SAW}

In this appendix, we establish rigorously the existence
of the infinite half-space SAW in all dimensions $d$.
In other words, we prove that the weak limit as $n\to\infty$ of
the uniform measure on $n$-step self avoiding walks in the half-space
$\{(x^1,\dots,x^d)\in \Z^d : x^1\ge 0\}$ starting from $0$ exists.

A nearest neighbor path $\omega=[\omega_0,\dots,\omega_n]$
in $\Z^d$ is called a {\em bridge} if
$\pi_1\omega_1<\pi_1\omega_j\le\pi_1\omega_n$ for all
$j=2,3,\dots,n$, where $\pi_1(x^1,\dots,x^d)=x^1$ denotes the projection
onto the first coordinate.  Madras and Slade~\cite{MS} 
used  Kesten's results  \cite {Kestensaw} to
 prove that
the uniform measure on $n$-step bridges starting at $0$ has a weak limit
as $n\to\infty$.  Our proof below is an adaptation of their methods.
In fact, the infinite half-space SAW is the same
as the limit measure for bridges.

We now introduce some notations and review some necessary background.
Let $\Upsilon_n$ denote the collection of $n$-step (nearest neighbor)
SAWs
 $\bigl[\omega_0,\dots,\omega_n]$
in $\Z^d$ with $\omega_0=0$ and satisfying
$\pi_1\omega_j>0$ for all $j=1,2,\dots,n$.
Set $\upsilon_n:=\#(\Upsilon_n)$, and
let $\mu_n$ denote the uniform measure on $\Upsilon_n$.
(If $d=2$, then $\upsilon_n = \#(\Lambda_n^+), $ where
$\Lambda_n^+$ is defined in \S\ref{halfsec}.)
A {\em renewal time} for $\omega\in\Upsilon_n$ is a $j\in \{1,2,\dots,n-1\}$
such that $\pi_1\omega_k\le\pi_1\omega_j<\pi_1\omega_m$ holds
for all $k$ and $m$ satisfying $0\le k\le j < m\le n$.
A bridge that has no renewal time is called {\em irreducible}.
Following the notations of~\cite{MS}, let $\lambda_n$ denote
the number of $n$-step irreducible bridges in $\Upsilon_n$.
A key step in establishing the limit for irreducible
bridges is Kesten's relation \cite {Kestensaw} (see also~\cite[(4.2.4)]{MS})
\begin{equation}\label{e.lambs}
\sum_{n=1}^\infty \lambda_n\beta^{-n}=1\,,
\end{equation}
where, as before, $\beta$ denotes the connective constant.
Kesten
also proved that 
\begin{equation}\label{e.2}
\lim_{n\to\infty}{\upsilon_{n+2}}/{\upsilon_n}=\beta^2\,;
\end{equation}
(see e.g., the proof of~\cite[Thm.~7.3.4]{MS}).
Shortly, we will justify the stronger statement
\begin{equation}\label{e.1}
\lim_{n\to\infty}{\upsilon_{n+1}}/{\upsilon_n}=\beta\,.
\end{equation}

For $\omega\in\Upsilon_n$, let $s(\omega)$ denote
the least renewal time of $\omega$.  If there is no
renewal time, then $s(\omega):=\infty$.
If $n>k$, then every walk  $\omega\in\Upsilon_n$ which satisfies
$s(\omega)=k$ is obtained by
concatenating an irreducible bridge of length $k$ with 
a (translated) element of $\Upsilon_{n-k}$. 
Conversely, every such concatenation
is an element of $\Upsilon_n$ satisfying $s(\omega)=k$.
Consequently,
$$
\mu_n\bigl\{s(\omega)=k\bigr\} =\lambda_k\, \upsilon_{n-k}/\upsilon_n\,.
$$
Hence,~(\ref{e.1}) gives
$$
\lim_{n\to\infty}\mu_n\bigl\{s(\omega)=k
\bigr\} =\lambda_k\, \beta^{-k}\,.
$$
Using~(\ref{e.lambs}) it now easily follows that the weak limit
$\lim_{n\to\infty} \mu_n$ is obtained by concatenating
an i.i.d.~sequence $\chi^1,\chi^2,\chi^3,\dots$ of
irreducible bridges, where the probability that $\chi^m$ is any particular
irreducible bridge of length $k$ is $\beta^{-k}$.

It remains to prove~(\ref{e.1}).
Note that
$$
\upsilon_n=\#(\Upsilon_n) \ge
\#\{\omega\in\Upsilon_n:s(\omega)\le k\}
=  \sum_{j=1}^k \lambda_j\upsilon_{n-j}\,.
$$
Therefore,
$$
1 \ge
\sum_{0<j<k/2}
\Bigl( \lambda_{2j}\frac{\upsilon_{n-2j}}{\upsilon_n}
+\lambda_{2j+1}\frac{\upsilon_{n-2j-1}}{\upsilon_{n-1}}
\frac{\upsilon_{n-1}}{\upsilon_n}\Bigr)\,.
$$
We take $\limsup$ of both sides as $n\to\infty$.  Using~(\ref{e.2})
this gives
$$
1\ge 
\sum_{0<j<k/2}
\Bigl( \lambda_{2j}\beta^{-2j}
+\lambda_{2j+1}
\beta^{-2j-1}
\limsup_{n\to\infty}\frac{\beta\upsilon_{n-1}}{\upsilon_n}\Bigr)\,.
$$
Taking $k$ to infinity and applying~(\ref{e.lambs}) now gives
$$
0\ge\Bigl( \sum_{j>0} \lambda_{2j+1} \beta^{-2j-1}\Bigr)
\limsup_{n\to\infty}\Bigl(\frac{\beta\upsilon_{n-1}}{\upsilon_n}-1\Bigr)\,.
$$
That is, $\limsup_{n\to\infty}( \upsilon_{n-1}/\upsilon_n) \le \beta^{-1}$.
This together with~(\ref{e.2}) implies~(\ref{e.1}), and completes the proof.

\medskip
\noindent {\bf Remark.}  If $a < \beta^{-1}$ and $\tilde \mu_a$ is
the probability
measure on $\Upsilon := \cup_n \Upsilon_n$
that assigns weight proportional to $ a^n$ to every $\omega
\in \Upsilon_n$, then   the existence of the
weak limit of $\tilde\mu_a$ as $a \nearrow \beta^{-1}$ can
be established
immediately from~(\ref{e.lambs}) without using~(\ref{e.2}) or~(\ref{e.1}).

\end{document}